\documentclass[11pt,twoside]{article}
\usepackage{amsmath,bm}
\usepackage{amssymb}
\usepackage{fancyhdr}
\usepackage{latexsym}
\usepackage{bbding}
\usepackage{mathrsfs}
\usepackage{exscale}
\usepackage{relsize}
\usepackage{wasysym}
\usepackage{cite}
\usepackage{multicol,graphics}
\usepackage{color}
\usepackage{graphicx}
\usepackage{subcaption}

\makeatletter
\@addtoreset{equation}{section}
\makeatother

\allowdisplaybreaks[4]
\def\to{\rightarrow}
\def\e{\left}
\def\r{\right}

\def\fr{\frac}

\def\i{\infty}

\def\R{\mathbb{R}}

\def\p{\partial}

\def\dl{\delta}
\def\vp{\varepsilon}

\def\bar{\overline}
\def\u{\underline}

\renewcommand\tilde{\widetilde}
\newcommand{\be}{\begin{equation}}
\newcommand{\ee}{\end{equation}}
\newcommand{\bc}{\begin{cases}}
\newcommand{\ec}{\end{cases}}
\newcommand{\bes}{\begin{equation*}}
\newcommand{\ees}{\end{equation*}}
\newcommand{\bls}{\begin{align*}}
\newcommand{\els}{\end{align*}}
\newcommand{\baa}{\begin{array}}
\newcommand{\eaa}{\end{array}}
\newcommand{\ba}{\begin{eqnarray}}
\newcommand{\ea}{\end{eqnarray}}
\newcommand{\bas}{\begin{eqnarray*}}
\newcommand{\eas}{\end{eqnarray*}}
\newcommand{\bd}{\begin{description}}
\newcommand{\ed}{\end{description}}

\newtheorem{theo}{\bf Theorem}[section]
\newtheorem{lem}[theo]{\bf Lemma}
\newtheorem{prop}[theo]{\bf Proposition}
\newtheorem{cor}[theo]{\bf Corollary}
\newtheorem{defi}[theo]{\bf Definition}
\newtheorem{rem}[theo]{\bf Remark}

\newenvironment{pr}[1][Proof]{\noindent\textbf{#1.} }{\hfill $\Box$}
\allowdisplaybreaks
\begin{document}
\date{}
\title{\bf{The speeds of propagation for the monostable Lotka-Volterra competition-diffusion system in general unbounded domains}}
\author{ Yang-Yang Yan$^{a,b}$ \  Wei-Jie Sheng$^{a,}$\thanks{Corresponding author
 (E-mail address: shengwj09@hit.edu.cn.)} \
\\
\footnotesize{$^a$ School of Mathematics, Harbin Institute of Technology}, \\
\footnotesize{Harbin, Heilongjiang, 150001, People's Republic of China}\\
\footnotesize{$^b$ Aix Marseille Univ, CNRS, I2M, Marseille, France}\\
}

\maketitle

\begin{abstract}
This paper is concerned with the speeds of propagation for the monostable Lotka-Volterra competition-diffusion system in general unbounded domains of $\mathbb{R}^N$.
We first establish various definitions of spreading speeds at large time in the situation where one species is an invader and the other is a resident.
Then, we study fundamental properties of these new definitions, including their relationships and their dependence on the  geometry of the domain and the initial values.
Under the conditions that both species possess the same diffusion ability and that the interactions between them are sufficiently weak,
we derive an upper bound for the spreading speeds in a large class of domains.
Furthermore, we obtain general upper and lower bounds for spreading speeds in exterior domains, as well as a general lower bound in domains containing large half-cylinders.
Finally, we construct some particular domains for which the spreading speeds can be zero or infinite.
\end{abstract}

\textbf{Keywords}: Spreading speeds; unbounded domains; monostable;  Lotka-Volterra 
competition-diffusion systems.

\textbf{AMS Classification}: 35B30; 35B40; 35K51; 35K57.

\section{Introduction and main results}
\subsection{Background and motivation}
In this paper, we focus on the speeds  of propagation for the following
Lotka-Volterra two-species competition-diffusion system
\begin{align}\label{clv}
\begin{cases}
  u_t =d_1\Delta u +r_1u (1-u -a_1v )&\text{ in }\ (0,+\i)\times\Omega, \\
v_t =d_2\Delta v +r_2v (1-v -a_2u )&\text{ in }\ (0,+\i)\times\Omega, \\
(u_\nu, v_\nu) =(0,0)&\text{ on }\ (0,+\i)\times\p\Omega, \\
(u,v)(0,\cdot)=(u_0,v_0)&\text{ in }\ \Omega,
\end{cases}
\end{align}
where  $u=u(t,x)$ and $v=v(t,x)$ represent the densities of two competitive species at time $t$ and location $x$, respectively.
The domain $\Omega$ is an open connected subset in $\R^N$ {($N\geq2$)} with smooth boundary $\p\Omega$,
and $\nu=\nu(x)$ is the outward unit normal for the point $x\in\p\Omega$.
The Neumann boundary condition $(u_\nu, v_\nu) =(0,0)$ implies that there is no flux cross $\p\Omega$.
The parameters $d_1$, $d_2$, $r_1$, $r_2$, $a_1$ and  $a_2$ are positive real numbers,
where $d_1$, $r_1$, $a_1$ and $d_2$, $r_2$, $a_2$ denote the diffusion rates, intrinsic growth rates, and competition coefficients of $u$ and $v$, respectively.

Throughout this paper, we always assume that the initial value $(u_0,v_0)$ satisfies the condition
\begin{align*}
u_0,v_0\in C(\bar\Omega,[0,1]),\ \ u_0\not\equiv0
\ \text{ and }\ (u_0,v_0)=(0,1)\text{ outside a compact subset in }\bar\Omega.
\end{align*} One calls $\Theta$ the set of such functions $(u_0,v_0)$.
In biology, $(u_0,v_0)\in\Theta$ means that the species $u$ is an invader, while
the species $v$ is a resident that occupies the entire territory $\Omega$.
Furthermore, we always provide the monostable assumption
\begin{align}\label{monostable}
0<a_1<1<a_2,
\end{align}
which implies that $(1, 0)$ is a stable equilibrium and $(0,1)$ is an unstable equilibrium  of the corresponding kinetic system of \eqref{clv}.
This assumption
further indicates that the species $u$ has a stronger competition ability than the
species $v$.
Under the above assumptions, it follows from \cite[Theorem 1]{MS} that
Cauchy problem \eqref{clv} admits a unique solution $(u,v)\in \e(C^{0,1}((0,+\i)\times\bar\Omega)\cap C^{1,2}((0,+\i)\times\Omega)\r)^2$
for given $(u_0,v_0)\in\Theta$.

In the case that $\Omega=\R^N$,
the planar traveling front
is a fundamental tool for characterizing the invasion of species $v$ by species $u$.
For the  competition-diffusion system
 \begin{align}\label{sys}
\begin{cases}
  u_t =d_1\Delta u +r_1u (1-u -a_1v )&\ \text{ in }\ \R\times\R^N, \\
v_t =d_2\Delta v +r_2v (1-v -a_2u )&\ \text{ in }\ \R\times\R^N,
\end{cases}
\end{align}
a planar traveling  front of \eqref{sys} connecting $(0,1)$ and $(1,0)$ is a solution $(u,v)$ with special form
 $(u,v)(t,x)=(\phi,\psi)(x\cdot e+c t)$,
 where $(\phi(\xi),\psi(\xi))$ ($\xi=x\cdot e+c t$)
 satisfying
 \begin{align}\label{ptr}
 \begin{cases}
d_1\phi''-c \phi'+r_1\phi(1-\phi-a_1\psi)=0,&\xi\in\R,\\
d_2\psi''- c \psi'+r_2\psi(1-\psi-a_2\phi)=0,&\xi\in\R,\\
 \phi'>0,\ \psi'<0, &\xi\in\R,\\
 (\phi,\psi)(-\i)=(0,1),\  (\phi,\psi)(+\i)=(1,0)
 \end{cases}
 \end{align}
 is the wave profile, $c\in\R$ is the wave speed and $-e\in\mathbb S^{N-1}$ (where $\mathbb S^{N-1}$ is the unit sphere of $\mathbb R^N$) is the propagation direction.
 It is well known from \cite{VVV} that there is a positive number $c_*$ such that
 system \eqref{sys} admits a planar traveling front of the form \eqref{ptr}
 if and only if $c\geq c_*$.
 The number $c_*$ is called the minimal wave speed of system \eqref{ptr}.
Geometrically, the level sets of such planar traveling fronts are parallel hyperplanes orthogonal to the
propagation direction $-e$, and the fronts are invariant in the moving frame with speed $c$.
 Biologically, the planar traveling fronts represent the spatial invasion process where species $u$ competitively excludes species $v$.
More precisely, the existence of such fronts  indicates that there is a transition zone
 moving from the steady state with the absence of species $u$
 to the steady state with the absence of species $v$.
 For more results about the planar traveling fronts of system \eqref{sys}, one can refer to
\cite{AX,CG,Ga,K1,KF,TF}  and references therein.

Regarding the minimal wave speed $c_*$, Kan-on \cite{K} showed that  $ 2\sqrt{d_1r_1(1-a_1)}\leq c_*\leq 2\sqrt{d_1r_1}$,
where $2\sqrt{d_1r_1(1-a_1)}$
 is the linear speed derived by linearizing \eqref{ptr} at the unstable state $(0,1)$.
In particular,  we say that $c_*$ is \textit{linearly determined} when  $c_*=2\sqrt{d_1r_1(1-a_1)}$
and $c_*$ is \textit{nonlinearly determined} when  $c_*>2\sqrt{d_1r_1(1-a_1)}$.
The speed determinacy mechanism of \eqref{sys} has been studied extensively in the past decades.
Hosono\cite{Ho} used numerical tests to verify that linear determinacy can hold or fail under proper conditions.
Later, Lewis, Li and Weinberger \cite{llw} proved that
$c_*$ is linearly determined under the conditions $0<d_2\leq2d_1$ and
\begin{align}\label{ls1}
(a_1,a_2,r_1,r_2)\in\left\{a_1a_2\leq1\right\}\cup\left\{a_1a_2>1,0<r_2\leq\frac{r_1(2-d_2/d_1)(1-a_1)}{a_1a_2-1}\right\}.
\end{align}
Huang \cite{Hu} indicated that the linear determinacy is realized under the condition
\begin{align*}
\frac{r_1(2-d_2/d_1)(1-a_1)+r_2}{r_2a_2}\geq\max\left(a_1,\frac{d_2-2d_1}{2|d_2-d_1|}\right).
\end{align*}
In addition, sufficient conditions on nonlinear determinacy have also received great attention.
Huang and Han \cite{HH} showed that $c_*$ is not linearly determined if $d_1r_2=d_2r_1$ and  $a_1\in[1-\vp,1)$ for some small $\vp>0$.
Alhasanat and Ou \cite{Al} found that $c_*$ is not linearly determined when
$$
\frac{r_1(d_2/d_1+2)(1-a_1)+r_2}{r_2a_2}<1-2(1-a_1),
$$
and $c_*$ is  linearly determined if $a_1<1/3$ and
$$
\frac{r_1(d_2/d_1-4)(1-a_1)}{4}<r_2<\frac{d_2r_1(1-a_1)}{2d_1a_2}
\ \text{ or }\ \frac{d_2r_1(1-a_1)}{2d_1a_2}\leq r_2<\frac{r_1(d_2/d_1+4)(1-a_1)}{4(a_2-1)}.
$$
It is worth noting that  a sharp necessary and sufficient condition on the
parameters that guarantees linear determinacy remains undiscovered.
For more discussion on linear/nonlinear determinacy, see \cite{AO,GL,HO,MO,RH} and references therein.

The asymptotic spreading speed also plays an important role in understanding  the large time behavior of solutions to system \eqref{clv}
when $\Omega=\R$.
For the following Cauchy problem
\begin{align}\label{clvr}
\begin{cases}
  u_t =d_1\Delta u +r_1u (1-u -a_1v ) & \ \text{ in }\ (0,+\i)\times\R, \\
v_t =d_2\Delta v +r_2v (1-v -a_2u ) & \ \text{ in }\ (0,+\i)\times\R, \\
(u,v)(0,\cdot)=(u_0,v_0) & \ \text{ in }\ \R,
\end{cases}
\end{align}
Weinberger, Lewis and Li \cite{wll} proved that if  $0\leq u_0<1$, $u_0\not\equiv0$, $0<v_0\leq 1$, $(u_0,v_0)=(0,1)$
 outside a compact subset in $\R$, then both species $u$ and $v$ spread at the
  same asymptotic spreading speed coinciding with the minimal wave speed $c_*$ in all directions for large time, that is,
\begin{align}\label{as}
\bc
\forall\ c>c_*,&
\sup\limits_{|x|\geq ct}\left(|u(t,x)|+|1-v(t,x)|\right)\to0,\\
\forall\ 0\leq c<c_*,&
\sup\limits_{|x|\leq ct}\left(|1-u(t,x)|+|v(t,x)|\right)\to0
\ec
\ \text{ as }{t\to+\i}.
\end{align}
Note that \eqref{as} means that an observer who  moves to left or right at a fixed speed greater than $c_*$ will see the extinction of species $u$ and the  persistence of species $v$, and that an observer who  moves to left or right at a fixed speed less than $c_*$ will see the extinction of species $v$ and the  persistence of species $u$.
In other words, the asymptotic speed $c_*$ describes the spreading speed of the invasion to the unstable equilibria $(0,1)$  by the stable equilibria $(1,0)$.
Girardin and Lam \cite{GLam} studied the spreading speeds of the solution to Cauchy problem \eqref{clvr}
when the initial value $(u_0,v_0)$ are null
or exponentially decaying  on the right half line.
By virtue of technical sub- and supersolutions,
they provided a complete understanding of spreading  properties of two invaders $u$ and $v$.
For more results on the spreading speeds of the solution to Cauchy problem \eqref{clvr},
see \cite{C,LL,LLL1,LLL2,PWZ} and references therein.

In the real world, the territories are often characterized by intricate and complex geometries,
which makes the study of the propagation phenomena in general unbounded domains vital and essential.
Note that although planar traveling  fronts no longer exist in general unbounded domains, the notions can be generalized to transition fronts in unbounded domains with arbitrary geometries, see \cite{BH1,BH2,shen}.
In the celebrate work \cite{bhn}, Berestycki, Hamel and Nadirashvili generalized the notions of asymptotic spreading speeds
for classical Fisher-KPP scalar equation
\begin{align}\label{FK}
\bc
u_t=\Delta u+f(u)&\ \text{ in }\ (0,+\i)\times\Omega,\\
u_\nu=0&\ \text{ in }\ (0,+\i)\times\p\Omega,\\
u(0,\cdot)=u_0&\ \text{ on }\ \Omega
\ec
\end{align}
 in general unbounded domain $\Omega$,
 where the nonlinearity $f\in C^1(\R_+)$ satisfies
\begin{align*}
 \bc
 f(0)=f(1)=0,\ f'(0)>0,\  f'(1)<0,\  f>0 \text{  in }(0,1),\
f<0\text{ in }(1,+\i),\\
f(s)\leq f'(0)s\  \text{ for all }s\in[0,1],
\ec
\end{align*}
 and the initial value $u_0$ is continuous, nonnegative, nonzero and  compactly supported in $\bar\Omega$.
More precisely,
they provided various new definitions of spreading speeds for solutions of \eqref{FK} at large time
and studied the relationships between different definitions.
They also analyzed the effect of the geometry of the domain and the initial value on these spreading speeds.
In addition, they proved that these new spreading speeds can be compared with the asymptotic spreading speed $2\sqrt{f'(0)}$
 in $\R^N$ under certain appropriate
geometry conditions.
It is worth pointing that the spreading speeds of Fisher-KPP equations with periodic coefficients in domains with periodic structures are
studied in  \cite{bhn1}.

In this paper, motivated by \cite{bhn},  we consider the spreading properties  of
 the monostable Lotka-Volterra two-species competition-diffusion system in general unbounded domains.
 More precisely, we first generalize the definitions of asymptotic spreading speeds for the solution $(u,v)$
 of Cauchy problem \eqref{clv} with initial value $(u_0,v_0)\in \Theta$ in $\bar\Omega$, under the assumption \eqref{monostable}.
 We then analyze  fundamental properties of these new definitions,
 such as relationships between different definitions,
 their dependence on  initial values, their upper and lower bounds and so on.
 In particular, when $\Omega$ is an exterior domain,
 we  derive general upper and lower bounds for the spreading speeds and
 establish a sufficient condition of parameters that ensures the spreading speeds coincide with
 the minimal wave speed $c_*$ of \eqref{ptr}.
 In addition, we  demonstrate that the spreading speed may be zero or infinite in certain special domains
 by providing suitable parameter conditions.

\subsection{ Spreading speeds in general unbounded domains}
Let us present the precise definitions of spreading speeds in general unbounded domains.
To this end, we first recall the so-called ``strongly unbounded'' domains.
In the sequel, the symbol $B(x,r)$ represents the open Euclidean ball of center $x\in\R^N$ and radius $r>0$.

\begin{defi}[\!\cite{bhn}]\label{unb}
  We say that a connected open set $\Omega\subset\mathbb R^N$ is strongly unbounded in a direction $e\in \mathbb S^{N-1}$
  if there exist $R_0\geq0$ and $s_0\in\mathbb R$ such that $\bar{B(se,R_0)}\cap\bar\Omega\neq\emptyset$ for all $s\geq s_0$. In particular, we define
  $$
  R(e)=\inf\{R\geq0:\exists\ s\in\R,\ \forall\ s'\geq s,\ \bar{B(s'e,R)}\cap\bar\Omega\neq\emptyset\}.
  $$
\end{defi}

Obviously, the whole space $\R^N$ and the exterior domain are two typical examples of being strongly
unbounded in all directions $e\in\mathbb S^{N-1}$.

  \begin{defi}\label{s1}
  Let $e\in\mathbb S^{N-1}$ be a direction in which $\Omega$ is strongly unbounded and let $R(e)\geq0$ be as in Definition \ref{unb}.
  Let $(u,v)$ be the solution of problem \eqref{clv} with initial value $(u_0,v_0)\in\Theta$.
  We define the upper spreading speed of $(u,v)$ in the direction $e$ as
  $$
  w^*(e,u_0,v_0)=\inf\left\{c>0:\forall\ A> R(e),\ \limsup_{t\to+\i}
  \left(\sup_{\substack{s\geq ct,\\ x\in\bar{B(se,A)}\cap\bar\Omega}}
  \left(|u(t,x)|+|1-v(t,x)|\right)\right)=0\right\}.
  $$
 We set $w^*(e,u_0,v_0)=+\i$ if there is no $c>0$ such that
$$
\sup_{s\geq ct,\ x\in\bar{B(se,A)}\cap\bar\Omega}
  \left(|u(t,x)|+|1-v(t,x)|\right)\to0\
  \text{ as }t\to+\i\text{ for all }A>R(e).
$$
\end{defi}

\begin{defi}\label{ss1}
Under the same conditions as in Definition \ref{s1},
assume that $(u,v)(t,x)\to(1,0)$ locally uniformly in $x\in\bar\Omega$ as $t\to+\i$.
  We define the lower spreading speed of $(u,v)$ in the direction $e$ as
  $$
  w_*(e,u_0,v_0)=\sup\left\{c>0:\forall\ A>R(e),\ \limsup_{\substack{\tau,t\to+\i,\\ \tau\leq ct}}\left(\sup_{\substack{\tau\leq s\leq ct,\\ x\in\bar{B(se,A)}\cap\bar\Omega}}
  \left(|1-u(t,x)|+|v(t,x)|\right)\right)=0\right\}.
  $$
We set $w_*(e,u_0,v_0)=0$ if there is no $c>0$ such that
$$
\limsup_{\tau,t\to+\i,\ \tau\leq ct}
\left(\sup_{\tau\leq s\leq ct,\ x\in\bar{B(se,A)}\cap\bar\Omega}
\left(|1-u(t,x)|+|v(t,x)|\right)\right)=0\text{ for all }A> R(e).
$$
  \end{defi}

The nonnegative real numbers $w^*(e,u_0,v_0)$ (if it is finite) and $w_*(e,u_0,v_0)$ can be viewed as, respectively,
the asymptotic spreading
speeds of the leading edge (where the reaction starts) and  the expanding region (where $(u,v)$ converges to $(1,0)$) for the solution $(u,v)$  uniformly with respect to all cylinders
along the direction $e$.
 When $w_*(e,u_0,v_0)=w^*(e,u_0,v_0)$, we say
that there exists an exact spreading speed of the solution $(u,v)$  uniformly in all cylinders along the direction  $e$.
Notice that the additional assumption in Definition \ref{ss1} is introduced to identify
the possible steady states.

\begin{defi}\label{s2}
 Under the same conditions as in Definition \ref{s1}, for $z\in\mathbb R^N$,
 we define the upper spreading speed of $(u,v)$ along the half-line $z+\mathbb R_+ e$
 as
  $$
  w^*(e,z,u_0,v_0)=\inf\left\{c>0:\exists\ A> 0,\ \limsup_{t\to+\i}\left(\sup_{\substack{s\geq ct,\\ x\in\bar{B(z+se,A)}\cap\bar\Omega}}
  \left(|u(t,x)|+|1-v(t,x)|\right)\right)=0\right\}.
  $$
We set $w^*(e,z,u_0,v_0)=+\i$ if for all $c>0$ and $A>0$, there holds
$$
\sup_{s\geq ct,\ x\in\bar{B(z+se,A)}\cap\bar\Omega}
  \left(|u(t,x)|+|1-v(t,x)|\right)\not\to0\text{ as }t\to+\i.
$$
\end{defi}

\begin{defi}\label{ss2}
Under the same conditions as in Definition \ref{ss1}, for $z\in\mathbb R^N$,
we define the lower spreading speed of $(u,v)$ along the half-line $z+\mathbb R_+ e$
 as
   $$
  w_*(e,z,u_0,v_0)=\sup\left\{c>0:\exists\ A> 0, \limsup_{\substack{\tau,t\to+\i,\\ \tau\leq ct}}
  \left(\sup_{\substack{\tau\leq s\leq ct,\\ x\in\bar{B(z+se,A)}\cap\bar\Omega}}
  \left(|1-u(t,x)|+|v(t,x)|\right)\right)=0\right\}.
  $$
We set
$w_*(e,z,u_0,v_0)=0$ if for all $c>0$ and $A>0$, there holds
$$
\limsup_{\tau,t\to+\i,\ \tau\leq ct}\left(\sup_{\tau\leq s\leq ct,\ x\in\bar{B(z+se,A)}\cap\bar\Omega}
  \left(|1-u(t,x)|+|v(t,x)|\right)\right)\neq0.
$$
  \end{defi}

  As previously mentioned,  the nonnegative real numbers $w^*(e,z,u_0,v_0)$
  (if it is finite) and $w_*(e,z,u_0,v_0)$ can be interpreted as the asymptotic spreading
speeds of the leading edge and the expanding region for the solution $(u,v)$
locally along the half-line $z+\mathbb R_+e$, respectively.
If $w_*(e,z,u_0,v_0)=w^*(e,z,u_0,v_0)$, we say that the solution
$(u,v)$ exhibits an exact spreading speed locally along the half-line $z+\mathbb R_+e$.

\begin{rem}\label{rks}\rm{
Denote
$$
R(e,z)=\inf\{R\geq0:\exists\ s\in\R,\ \forall\ s'\geq s,\ \bar{B(z+s'e,R)}\cap\bar\Omega\neq\emptyset\}.
$$
It is easy to see that $|R(e)-|z-(z\cdot e) e||\leq R(e,z)\leq R(e)+|z-(z\cdot e) e|$
for all $z\in\R^N$ (see Figure \ref{fig1-2})
\begin{figure}[htbp]
    \centering
    \begin{subfigure}[b]{0.49\textwidth}
        \centering
        \includegraphics[width=7.8cm]{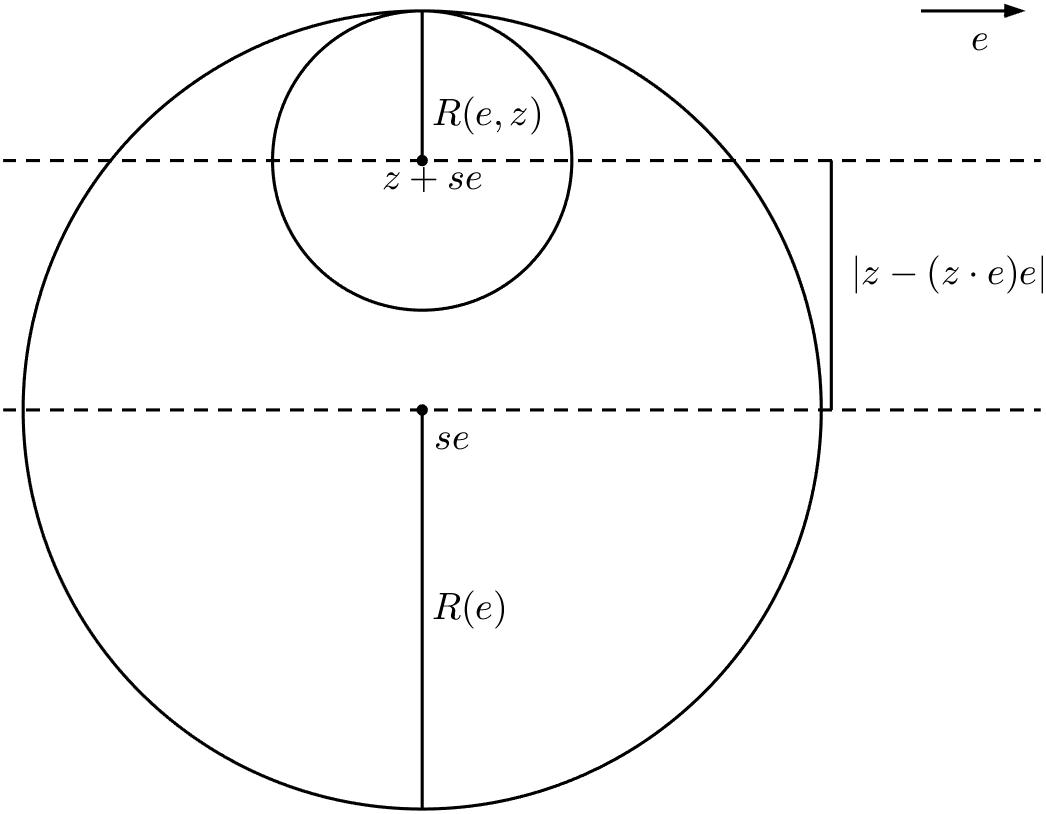}
    \end{subfigure}
    \hfill
    \begin{subfigure}[b]{0.49\textwidth}
        \centering
        \includegraphics[width=8cm]{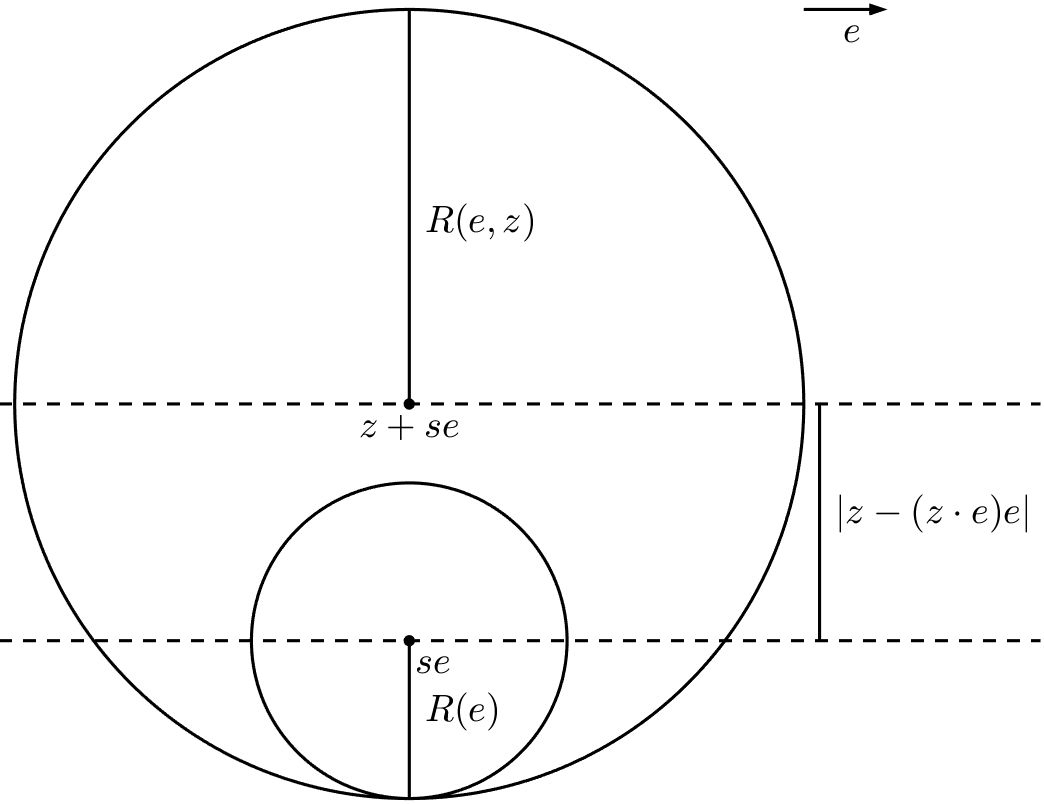}
    \end{subfigure}
    \caption{the relationship between $R(e)$ and $R(e,z)$.}\label{fig1-2}
\end{figure}
 and
 $R(e,z)=R(e)$ if $z$ is the origin in $\R^N$.
  If $R(e,z)>0$ and there exists $s\in\R$ such that
$\bar{B(z+s'e,R(e,z))}\cap\bar\Omega\neq\emptyset$ for all $s'\geq s$, then  $$
  w^*(e,z,u_0,v_0)=\inf\left\{c>0:\ \limsup_{t\to+\i}\left(\sup_{\substack{s\geq ct,\\ x\in\bar{B(z+se,R(e,z))}\cap\bar\Omega}}
  \left(|u(t,x)|+|1-v(t,x)|\right)\right)=0\right\}
  $$
 and
   $$
  w_*(e,z,u_0,v_0)=\sup\left\{c>0:\ \limsup_{\substack{\tau,t\to+\i,\\
  \tau\leq ct}}\left(\sup_{\substack{\tau\leq s\leq ct,\\ x\in\bar{B(z+se,R(e,z)}\cap\bar\Omega}}
  \left((|1-u(t,x)|+|v(t,x)|\right)\right)=0\right\}.
  $$
 On the other hand,  if $R(e,z)=0$ or there is no $s\in\R$ such that
$\bar{B(z+s'e,R(e,z))}\cap\bar\Omega\neq\emptyset$ for all $s'\geq s$ (see \cite[Remark 1.4]{bhn}), then
\begin{align*}
  w^*(e,z,u_0,v_0)=\inf\left\{c>0:\exists A>R(e,z), \limsup_{t\to+\i}\left(\sup_{\substack{s\geq ct,\\ x\in\bar{B(z+se,A)}\cap\bar\Omega}}
  \left(|u(t,x)|+|1-v(t,x)|\right)\right)=0\right\}
\end{align*}
 and
   $$
  w_*(e,z,u_0,v_0)=\sup\left\{c>0:\exists A>R(e,z),\limsup_{\substack{\tau,t\to+\i,\\ \tau<ct}}
  \left(\sup_{\substack{\tau\leq s\leq ct,\\ x\in\bar{B(z+se,A)}\cap\bar\Omega}}
  \left(|1-u(t,x)|+|v(t,x)|\right)\right)=0\right\}.
  $$
Moreover, it follows from Definitions \ref{s1}-\ref{ss2} that
  \begin{align*}
  &\forall\ \gamma> w^*(e,u_0,v_0),\ \forall\ A> R(e), \ \sup_{x\in\bar{ B(\gamma t e,A)}\cap\bar\Omega}(|u(t,x)|+|1-v(t,x)|)\to0\ \text{ as }t\to+\i,\\
  &\forall\ 0<\gamma< w_*(e,u_0,v_0),\ \forall\ A> R(e),\ \sup_{x\in\bar{ B(\gamma t e,A)}\cap\bar\Omega}(|1-u(t,x)|+|v(t,x)|)\to0\ \text{ as }t\to+\i,\\
  &\forall\ \gamma> w^*(e,z,u_0,v_0),\ \exists\ A>0,\ \sup_{x\in\bar{ B(\gamma t e,A)}\cap\bar\Omega}(|u(t,x)|+|1-v(t,x)|)\to0\ \text{ as }t\to+\i,\\
  &\forall\ 0<\gamma< w_*(e,z,u_0,v_0),\ \exists\ A>0,\ \sup_{x\in\bar{ B(\gamma t e,A)}\cap\bar\Omega}(|1-u(t,x)|+|v(t,x)|)\to0\ \text{ as }t\to+\i.
  \end{align*}
}\end{rem}

\subsection{Fundamental properties}\label{fund}
We first focus on the relationship between these new spreading speeds
and their dependence on the position $z$.
Notice that it immediately follows from  Definitions \ref{s1}-\ref{ss2} that
\begin{align}\label{w}
0\leq w_*(e,u_0,v_0)\leq w_*(e,z,u_0,v_0)\leq w^*(e,z,u_0,v_0)\leq w^*(e,u_0,v_0)
\end{align}
for all $z\in\R$ and $e\in\mathbb S^{N-1}$.
In the following result, we present a more precise relationship between these spreading speeds when $z$ varies
and show a strict inequality in some special domains.

\begin{prop}\label{dz}
Let $N\geq2$ and $e\in\mathbb S^{N-1}$ be given.
Then the following statements hold:
\begin{itemize}
  \item [\rm(a)]For each locally $C^2$ domain $\Omega$
which is strongly unbounded in the direction $e$ and for each initial value $(u_0,v_0)\in\Theta$, one has
\begin{align}\label{dzz}
\sup_{z\in\R^N}w^*(e,z,u_0,v_0)=w^*(e,u_0,v_0).
\end{align}
In particular, given $z\in\R^N$, if $d_1= d_2$ and $a_1a_2\leq1$, then there are some locally $C^2$ domains $\Omega$ which are strongly unbounded in the direction $e$
 such that
  $w^*(e,z,u_0,v_0)<w^*(e,u_0,v_0)$
  for all $(u_0,v_0)\in\Theta$.
  \item [\rm(b)]Further assume that $(u,v)(t,x)\to(1,0)$ locally uniformly in $x\in\bar\Omega$ as $t\to+\i$.
  Then for each locally $C^2$ domain $\Omega$
which is strongly unbounded in the direction $e$ and for each initial value $(u_0,v_0)\in\Theta$, one has
\begin{align}\label{dzzz}
\inf_{z\in\R^N}w_*(e,z,u_0,v_0)=w_*(e,u_0,v_0).
\end{align}
\end{itemize}

\end{prop}

In fact, we prove that in the domain whose complement has the shape of a comb with larger and larger teeth,
if $d_1= d_2$ and $a_1a_2\leq1$, then
   $w^*(e,z,u_0,v_0)=0$ for some $z\in\R^N$   but $w^*(e,u_0,v_0)\geq2\sqrt{d_1r_1(1-a_1)}$.
   The parameter conditions $d_1= d_2$ and $a_1a_2\leq1$ emerge from the application of Theorem \ref{upd} in proving the second assertion of part (a).

   However, the speeds $w^*(e,z,u_0,v_0)$ and $w_*(e,z,u_0,v_0)$ does not always depend on the position $z$.
Under some  related connected assumptions, we provide a sufficient condition to ensure that
$w^*(e,z,u_0,v_0)$ and $w_*(e,z,u_0,v_0)$ are independent of the choice of $z$.
Before stating our results, let us first introduce the hypothesis $H_{y,z}$
and the interior ball property of $\Omega$.

{\bf Hypothesis $H_{y,z}$.} Let $\Omega$ be a strongly unbounded domain in a direction $e\in\mathbb S^{N-1}$ and $d_\Omega$ be the geodesic distance in $\bar\Omega$.
We say that the points $y$ and $z$ are asymptotically connected in the direction $e$ (or that hypothesis $H_{y,z}$ is satisfied) within $\Omega$
if there exist $R_y>R(e,y)$ and $R_z>R(e,z)$ such that
\begin{align*}
\limsup_{s\to+\i}\sup_{\substack{y'\in\bar{B(y+se,R_y)}\cap\bar\Omega\\
z'\in\bar{B(z+se,R_z)}\cap\bar\Omega}}d_\Omega(y',z')<+\i.
\end{align*}

{\bf Interior ball property.} We say that a domain $\Omega$ has the  interior ball property at point $x\in\p\Omega$, if there is an open
Euclidean ball $B\subset\Omega$ such that $x\in\p B$.
If $\Omega$ has the interior ball property at each point $x\in\p\Omega$, then we say that $\Omega$  has the interior ball property.
It is worth noting that the globally $C^2$ domain automatically satisfies the interior ball property.
Furthermore, if there exists $\vp_0>0$ such that for every point $x\in\p\Omega$, there is an open
Euclidean ball $B\subset\Omega$ with radius $\vp_0$ such that $x\in\p B$,
 then we say that $\Omega$ has the  interior ball property with a radius $\vp_0$.

\begin{theo}\label{idz}
 Let $N\geq2$ and $e\in\mathbb S^{N-1}$ be given. Let $\Omega$ be a strongly unbounded domain in the direction $e$.
Assume that there exists $\vp_0\in(0,1)$ such that $\Omega$ satisfies the interior ball property with a radius $\vp_0$,
and that the hypothesis $H_{y,z}$ is satisfied for some points $y$ and $z$ in $\R^N$.
  Then it holds
  \begin{align}\label{idzz}
\forall\ (u_0,v_0)\in\Theta,\ \ w^*(e,y,u_0,v_0)=w^*(e,z,u_0,v_0).
 \end{align}
Further assume that $(u,v)(t,x)\to(1,0)$ locally uniformly in $x\in\bar\Omega$ as $t\to+\i$,
then
  \begin{align}\label{idzzz}
\forall\ (u_0,v_0)\in\Theta,\ \   w_*(e,y,u_0,v_0)=w_*(e,z,u_0,v_0).
 \end{align}
 In particular, if the hypothesis $H_{y,z}$ is satisfied for all $y,z\in\R^N$, then \eqref{idzz} and \eqref{idzzz}
hold for all $y,z\in\mathbb R^N$ and $(u_0,v_0)\in\Theta$.
\end{theo}

Compared with the results for single Fisher-KPP equation,
we change the regularity assumption on the domain $\Omega$
from being globally $C^{2,\alpha}$ for some $\alpha > 0$ to requiring that $\Omega$ satisfies the interior ball property with a radius $\vp_0$.
The change arises from the absence of Harnack inequality at the boundary for parabolic differential systems.
In fact, this assumption implies that the domain is smooth in a certain sense and does not contain infinitely narrow channels.
Obviously, the domain with the shape of an infinite cusp (see the proof of Theorem \ref{00}) does not satisfy this assumption.

Next, let us study the dependence of the spreading speeds on  initial values.
By utilizing the comparison principle and the weak Harnack inequality,
we prove that the spreading speeds do not depend on the choice of initial values
under some appropriate assumptions.
\begin{theo}\label{idi}
Let $\Omega$ be a connected domain satisfying the interior ball property with a radius $\vp_0$ for some $\vp_0\in(0,1)$.
For each $z\in\mathbb R^N$ and $r>0$, let $\mu_r^z$ denote the Lebesgue measure of $\Omega\cap B(z,r)$.
Assume that there exists $R_0>0$ such that $\mu_r^z>0$ for all
$z\in\R^N$ and $r\geq R_0$, and that $\mu_{r+1}^z/\mu_r^z\to1$ as $r\to+\i$ uniformly in $z\in\mathbb R^N$.
Let $(u,v)$  be the solution of \eqref{clv} with a given initial value $(u_0,v_0)\in\Theta$.
Then
$$
(u,v)(t,x)\to(1,0)\ \text{ locally uniformly in }x\in\bar\Omega\text{ as }t\to+\i.
$$
Furthermore, $\Omega$ is strongly unbounded in any direction $e\in\mathbb S^{N-1}$.
If the initial value $(u_0,v_0)\in\Theta$ satisfies $u_0<1$ and $v_0>0$,
then $w^*(e,u_0,v_0)$, $w_*(e,u_0,v_0)$, $w^*(e,z,u_0,v_0)$ and $w_*(e,z,u_0,v_0)$ do not depend on
the choice of $(u_0,v_0)$.
\end{theo}

\begin{rem}{\rm
Since the weak Harnack inequality is only applied to Euclidean balls contained in $\Omega$,
the strong connectedness assumption (for all $x,y\in\bar\Omega$, if  $|x-y|<+\i$, then $d_\Omega(x,y)<+\i$) that made in \cite{bhn} is no longer necessary. Note that the whole space $\mathbb{R}^N$ and exterior domains of class $C^2$ satisfy the assumptions of Theorem \ref{idi}.
A further example is the domain $\Omega := \mathbb{R}^N \setminus \{x \in B(se, R_0) : s \geq 0\}$ for some $e \in \mathbb{S}^{N-1}$ and $R_0 > 0$, which satisfies the assumptions of Theorem \ref{idi} but is not strongly connected.}
\end{rem}

Another key issue involves estimates of the spreading speeds.
With the help of some precise heat kernel estimates, we establish a general upper bound for the spreading speeds on a large class of domains
satisfying the extension property under certain parameter conditions.
Before stating our results, we first introduce
the extension property of $\Omega$  (see \cite{D}).

{\bf Extension property.} We say that a non-empty open subset $\Omega$ of $\R^N$  satisfies the  extension property if,
 for all $1\leq p\leq+\i$, there exists a bounded linear map $E:W^{1,p}(\Omega)\to W^{1,p}(\R^N)$
such that $E(f)$ is an extension of $f$ from $\Omega$ to $\R^N$ for all $f\in W^{1,p}(\Omega)$.
A fundamental result given in \cite{st} states that a domain $\Omega$ has the extension property if $\p\Omega$
is smooth in the following sense,
namely,
there exist constants $\vp>0$ and $M>0$, an integer $k$ and a countable sequence of open sets $U_n$ such that:

(a) if $x\in\p\Omega$, then the Euclidean  ball $B(x,\vp)$ is contained in $U_n$ for some $n$;

(b) no point in $\R^N$ is contained in more than $k$ distinct sets $U_n$;

(c) for each $n$ there exists an isometry $T_n:\R^N\to\R^N$ and a function $\phi_n:\R^{N-1}\to\R$
such that $|\phi_n(y)-\phi_n(y')|\leq M|y-y'|$ for all $y,y'\in\R^{N-1}$. Moreover, $U_n\cap\Omega=U_n\cap T_n\Omega_n$,
where
  $$
  \Omega_n=\{(z_1,\cdots,z_N)\in\R^N: \phi_n(z_1,\cdots,z_{N-1})<z_N\}.
  $$
Therefore, any bounded domain and any  exterior domain
with smooth boundary has the extension property.

\begin{theo}\label{upd}
Let $\Omega$ be a locally $C^2$ connected open subset of $\R^N$ satisfying the extension property.
Assume that $\Omega$ is strongly unbounded in a direction $e$.
 Let $(u,v)$ be the solution of \eqref{clv} with  initial value $(u_0,v_0)\in\Theta$.
 The following statements are true:
 \begin{itemize}
   \item [\rm(a)] Assume that $(u,v)(t,x)\to(1,0)$ locally uniformly in $x\in\bar\Omega$ as $t\to+\i$. Then
  \begin{align}\label{updd}
   w_*(e,u_0,v_0)\leq w_*(e,z,u_0,v_0)\leq 2\sqrt{d_1r_1}
  \end{align}
  for all $z\in\R^N$ and
  \begin{align}\label{updu}
  \forall\ c> 2\sqrt{d_1r_1},\
  \sup_{|x|\geq ct,\ x\in\bar\Omega}u(t,x)\to0\ \text{ as }t\to+\i.
\end{align}

   \item [\rm(b)]If $d_1=d_2$ and $a_1a_2\leq1$, then
\begin{align}\label{upd1}
w^*(e,z,u_0,v_0)\leq w^*(e,u_0,v_0)\leq 2\sqrt{d_1r_1(1-a_1)}
\end{align}
for all $z\in \R^N$ and
\begin{align}\label{upd2}
  \forall \ c> 2\sqrt{d_1r_1(1-a_1)},\
  \sup_{|x|\geq ct,\ x\in\bar\Omega}(|u(t,x)|+|1-v(t,x)|)\to0\ \text{ as }t\to+\i.
\end{align}
 \end{itemize}

\end{theo}

In part (a), the formula \eqref{updu} is derived through a direct comparison with the results of single Fisher-KPP equation.
This implies that, at large times, the invader $u$ becomes extinct in the area where  $|x|/t>2\sqrt{d_1r_1}$
($|\cdot|$ denotes the Euclidean norm in $\R^N$).
Since \eqref{updu} is stronger than \eqref{updd} (see the proof of Theorem \ref{upd}),
 the minimal wave speed $2\sqrt{d_1r_1}$ of \eqref{ptr} with $\psi\equiv0$,
turns out to be a general upper bound of $w_*(e,z,u_0,v_0)$
in the domain satisfying the extension property.
In part (b), \eqref{upd2}  means that
 if the invader $u$ and the resident $v$ have the same  dispersal abilities ($d_1=d_2$),
and the interactions between them are sufficiently weak ($a_1a_2\leq1$),
then at large times, the species $u$ goes extinct while the species $v$ persists in the area where
 $|x|/t>2\sqrt{d_1r_1(1-a_1)}$.
Notice that \eqref{upd2} is stronger than \eqref{upd1}.
When $d_1=d_2$ and $a_1a_2\leq 1$, since \eqref{ls1} yields that the minimal wave speed $c_*$ is linearly determined,
it follows that  the upper spreading speeds in
 domains satisfying the extension property  are bounded from above by the minimal wave speed $c_*$.
However,
due to the uncertainty in the geometry of domain and the coupling between $u$ and $v$,
we are unable to derive a general upper bound of $w^*(e,u_0,v_0)$.
We leave it for future studies.

\subsection{Exterior domains and domains containing large half-cylinders}
In this section,  we study the bounds of the spreading speeds
in two types of unbounded domains: exterior domains and domains containing large half-cylinders.
An exterior domain  is referred to
a non-empty connected open subset $\Omega$ of $\R^N$ such that
$\R^N\backslash \Omega$ is compact.
The whole space $\R^N$ can be viewed as a special case of an exterior domain.
Applying sub- and supersolution method, we derive a general upper bound of the spreading speeds in  exterior domains.
\begin{theo}\label{upd-ex}
Let $\Omega$ be an exterior domain of class $C^2$.
Then
\begin{align}\label{ue1}
w_*(e,u_0,v_0)\leq w_*(e,z,u_0,v_0)\leq w^*(e,z,u_0,v_0)\leq w^*(e,u_0,v_0)\leq 2\sqrt{r_1}\max(\sqrt{d_1},\sqrt{d_2/2})
\end{align}
for all $z\in \R^N$, $(u_0,v_0)\in\Theta$ and $e\in\mathbb S^{N-1}$.
\end{theo}

On the other hand, by comparing with the case where $v\equiv1$, we obtain
a general lower bound of the spreading speeds in exterior domains.

\begin{theo}\label{ext}
Let $\Omega$ be an exterior domain of class $C^2$.
If $(u,v)$ solves \eqref{clv} with initial value $(u_0,v_0)\in\Theta$, then
\begin{align}\label{ext2}
\forall\ 0\leq c< 2\sqrt{d_1r_1(1-a_1)}, \ \ \sup\limits_{|x|\leq ct,\ x\in\bar\Omega}
(|1-u(t,x)|+|v(t,x)|)\to0
\ \ \text{ as }t\to+\i.
\end{align}
In particular, for all $e\in\mathbb S^{N-1}$,
$ z\in\R^N$ and $(u_0,v_0)\in\Theta$,
\begin{align}\label{ext1}
w^*(e,u_0,v_0)\geq w^*(e,z,u_0,v_0)\geq
w_*(e,z,u_0,v_0)\geq w_*(e,u_0,v_0)\geq2\sqrt{d_1r_1(1-a_1)}.
\end{align}
\end{theo}

The formula \eqref{ext2} implies that at large times, an observer moving with a  mean speed less than $2\sqrt{d_1r_1(1-a_1)}$ will see the persistence of invader $u$ and the extinction of resident $v$.
Since \eqref{ext2} is stronger than \eqref{ext1},
it follows that the linear speed $2\sqrt{d_1r_1(1-a_1)}$ serves as
a general lower bound of spreading speeds in exterior domains.

Notice that exterior domains with $C^2$ boundary are typical examples of domains satisfying the extension property.
We can partially extend the classical work of Weinberger, Lewis and Li \cite{wll}
to  exterior domains of class $C^2$.
Specifically, if $d_1=d_2$ and $a_1a_2\leq1$,  then \eqref{ls1}, \eqref{w}, Theorems \ref{upd} and \ref{ext}  yield that
 the spreading speeds on $C^2$ exterior domains coincide with the minimal wave speed
$c_*=2\sqrt{d_1r_1(1-a_1)}$ of \eqref{ptr}.
This is indeed shown in the following corollary.

\begin{cor}\label{cor1}
Let $\Omega$ be an exterior domain of class $C^2$.
If $d_1=d_2$ and $a_1a_2\leq1$, then for all $e\in\mathbb S^{N-1}$,
$ z\in\R^N$ and $(u_0,v_0)\in\Theta$,
$$
w_*(e,u_0,v_0)= w_*(e,z,u_0,v_0)= w^*(e,z,u_0,v_0)= w^*(e,u_0,v_0)=2\sqrt{d_1r_1(1-a_1)}.
$$
\end{cor}

We now turn to consider domains containing large half-cylinders.
For the domain which contains a semi-infinite cylinder in the direction $e$ with large section,
the upper spreading speeds $w^*(e,u_0,v_0)$ and $w^*(e,z,u_0,v_0)$ are bounded from below by a constant close to the linear speed $2\sqrt{d_1r_1(1-a_1)}$,
as stated in the following theorem.
However, due to the uncertainty of the boundary, we are unable to derive general upper and lower bounds for the spreading speeds by using similar arguments as those applied in exterior domains.

\begin{prop}\label{half-cyl}
For any $\vp\in(0,2\sqrt{d_1r_1(1-a_1)}]$,
there exists $R_0=R_0(\vp)>0$ such that if
\begin{align}\label{cyc}
\Omega\supset \mathcal C_{e,A,x_0,R}
:=\ &\{x\in\R^N: x\cdot e\geq A,\ |(x-x_0)-((x-x_0)\cdot e)e|<R\}
\end{align}
for some $e\in\mathbb S^{N-1}$, $A\in\R$ and $R>R_0$,
then for all $(u_0,v_0)\in \Theta$ and $z\in\R^N$ such that
$|z-x_0-((z-x_0)\cdot e)e|<R$,
\begin{align*}
w^*(e,u_0,v_0)\geq
w^*(e,z,u_0,v_0)\geq 2\sqrt{d_1r_1(1-a_1)}-\vp.
\end{align*}
In particular, if $\Omega$ contains a sequence of semi-infinite cylinders
$\e(\mathcal C_{e,A_n,x_{0,n},R_n}\r)_{n\in\mathbb N}$ with
$A_n\in \R$, $x_{0,n}\in\R^N$ and $R_n\to+\i$
as $n\to+\i$,
then $w^*(e,u_0,v_0)\geq 2\sqrt{d_1r_1(1-a_1)}$ for all $(u_0,v_0)\in\Theta$.
\end{prop}

A domain $\Omega$ satisfying \eqref{cyc} is that it
 contains a ``quarter of space'', that is,
\begin{align}\label{quarter}
\Omega\supset \{x\in\R^N: x\cdot e>A,\ x\cdot e'>B\}
\end{align}
 for some $(A,B)\in\R^2$ and $e,e'\in\mathbb S^{N-1}$ with $e\cdot e'=0$.
It follows from proposition \ref{half-cyl} that for all $\vp\in(0,2\sqrt{d_1r_1(1-a_1)}]$,
 there exists $R_0>0$ such that
 $ w^*(e,z,u_0,v_0)\geq2\sqrt{d_1r_1(1-a_1)}-\vp$
 for all $(u_0,v_0)\in\Theta$ and
 $z\in\bigcup_{x_0\in\R^N,\ R\geq R_0,\ x_0\cdot e'>B+R}B(x_0,R)$.
Consequently,
 $$
 w^*(e,u_0,v_0)\geq w^*(e,z,u_0,v_0)\geq2\sqrt{d_1r_1(1-a_1)}
 $$
 for all $(u_0,v_0)\in\Theta$ and $z\in\R^N$ such that $z\cdot e'>B$.
If $\Omega$ further satisfies the extension property,
then Theorem \ref{upd}  implies  that the upper spreading speeds $w^*(e,u_0,v_0)$ and $w^*(e,z,u_0,v_0)$ coincide with the minimal wave speed
$c_*=2\sqrt{d_1r_1(1-a_1)}$ of \eqref{ptr} when $d_1=d_2$ and $a_1a_2\leq1$, as stated in the following corollary.

\begin{cor}\label{coro2}
Assume that $\Omega$ be a locally $C^2$ connected open subset of $\R^N$ satisfying \eqref{quarter} and
the extension property.
If $d_1=d_2$ and $a_1a_2\leq1$,
 then $w^*(e,z,u_0,v_0)=w^*(e,u_0,v_0)=2\sqrt{d_1r_1(1-a_1)}$ for all $(u_0,v_0)\in\Theta$ and $z\in\R^N$ such that $z\cdot e'>B$.
 \end{cor}

\subsection{Domains with zero or infinite spreading speeds}
In the previous conclusions,
we studied several types of domains with positive and finite spreading speeds.
However, there are also some unbounded domains where the spreading speeds are zero or infinite.
\begin{theo}\label{00}
The following statements hold:
\begin{itemize}
  \item [\rm(a)]
 There are some locally $C^2$ domains of $\R^2$ which satisfy the extension property and are
  strongly unbounded in all directions $e\in\mathbb S^{1}$, such that if $d_1=d_2$ and $a_1a_2\leq1$, then
  $$
w^*(e,z,u_0,v_0)=w^*(e,u_0,v_0)=0
  $$
  for all $e\in \mathbb S^1$, $ z\in\R^2$ and $(u_0,v_0)\in\Theta$.
  \item [\rm(b)]
   There are some locally $C^2$ domains of $\R^N$ $(N\geq2)$ which are
  strongly unbounded in all directions $e\in\mathbb S^{N-1}$ but do not satisfy the extension property such that
  $$
w_*(e,u_0,v_0)=w_*(e,z,u_0,v_0)=w^*(e,z,u_0,v_0)=w^*(e,u_0,v_0)=+\i
  $$
  for all $e\in \mathbb S^{N-1}$, $ z\in\R^N$ and $(u_0,v_0)\in\Theta.$
\end{itemize}
\end{theo}

Part (a) implies that even if the domain satisfies the extension property,
the upper spreading speeds may be zero when $d_1=d_2$ and $a_1a_1\leq1$.
In the proof of part (a), we construct a spiral-shaped domain that satisfies the extension property.
For this domain,  the upper spreading speeds are zero in any direction when $d_1=d_2$ and $a_1a_1\leq1$.
Furthermore, Definitions \ref{s1} and \ref{s2} yield that the invasion of species $v$ by species $u$ will eventually fail in this domain if $d_1=d_2$ and $a_1a_1\leq1$.

Compared part (b) with Theorem \ref{upd},
one observes that
the finite upper bounds no longer exist without extension property.
In the proof of part (b), we show that a domain with an infinite cusp lacking the extension property
has the infinite spreading speeds.
Definitions \ref{s1}-\ref{ss2} imply that at large times, the invasion of species $v$ by species $u$ will success in this domain.

 {\bf Outline of the paper.}
 This paper is organized as follows. In Section \ref{s11}, we introduce some preliminaries.
 In Section \ref{S1}, we study fundamental properties of the spreading speeds, that is,
 we prove formulas \eqref{dzz} and \eqref{dzzz} in Proposition \ref{dz}, Theorems \ref{idz}, \ref{idi} and \ref{upd}.
Section \ref{S2} is devoted to the bounds of the spreading speeds
in exterior domains and domains containing large half-cylinders, namely,
we prove Theorems \ref{upd-ex} and \ref{ext}, and Proposition \ref{half-cyl}.
In Section \ref{S3}, we focus on  domains with zero or infinite spreading speeds,
that is, we prove the second assertion of part (a) in Proposition \ref{dz} and Theorem \ref{00}.

\section{Preliminaries}\label{s11}
In this section, we introduce some comparison principles and  weak Harnack inequalities.
For any vectors $(p_1,p_2)$ and $(q_1,q_2)$ in $\R^2$, the symbol $(p_1,p_2)\ll (q_1,q_2)$ means $p_i<q_i$ for each $i=1,2$,
and $(p_1,p_2)\leq (q_1,q_2)$ means $p_i\leq q_i$ for each $i=1,2$.
For the convenience of using the comparison principle, let us transform the competition-diffusion system \eqref{clv}
into a cooperation-diffusion system.
More precisely, by changing variables $u^*=u$, $v^*=1-v$, $u^*_0=u_0$ and $v^*_0=1-v_0$,
and dropping the symbol $*$ for convenience,
we obtain that system \eqref{clv} is equivalent to
\begin{align}\label{colv}
\begin{cases}
  u_t=d_1\Delta u+r_1u(1-a_1-u+a_1v) &\text{ in }\ (0,+\i)\times\Omega, \\
v_t=d_2\Delta v+r_2(1-v)(a_2u-v) &\text{ in }\ (0,+\i)\times\Omega, \\
(u_\nu, v_\nu)=(0,0) &\text{ on }\ (0,+\i)\times\p\Omega, \\
(u,v)(0,\cdot)=(u_0, v_0)&\text{ in }\ \Omega.
\end{cases}
\end{align}
The constant equilibria $(0,1)$ and $(1,0)$
become $(0,0)$ and $(1,1)$, respectively. Clearly, $(0,0)$ is unstable equilibrium,
while $(1,1)$ is stable equilibrium of the corresponding kinetic system of \eqref{colv}.
According to the definition of set $\Theta$, one infers that the initial value $(u_0,v_0)$ satisfies
$$
u_0,v_0\in C(\bar\Omega,[0,1]),\ u_0\not\equiv0\text{ and } (u_0,v_0) \text{ is compact support in }\bar\Omega.
$$
Applying the maximum principle and Hopf boundary lemma, one concludes that
the solution $(u,v)$  of system \eqref{colv} satisfying
\begin{align}\label{0uv1}
(0,0)\leq(u,v)\leq(1,1) \text{ in }[0,+\i)\times\bar\Omega\ \text{ and }\
(0,0)\ll (u,v)\ll(1,1)\text{ in } (0,+\i)\times\bar\Omega.
\end{align}

We now state the definitions of sub- and supersolutions of system \eqref{colv}.
\begin{defi}\label{dss}
If a vector-value function $( u,v)$ satisfies
$ u, v\in C^{0,1}((0,+\i)\times\bar\Omega)\cap C^{1,2}((0,+\i)\times\Omega)$ and
such that
$$
\bc
\mathcal L_1(u,v):=  u_t-d_1\Delta u-r_1u(1-a_1-u+a_1v)\geq 0\ ( \leq0) &\text{ in }\ (0,+\i)\times\Omega, \\
\mathcal L_2(u,v):= v_t-d_2\Delta v-r_2(1-v)(a_2u-v)\geq0\ (  \leq0) &\text{ in }\ (0,+\i)\times\Omega,\\
(u_\nu,v_\nu) \geq(0,0)\ (  \leq(0,0)) &\text{ on }\ (0,+\i)\times\p\Omega,
\ec
$$
then $( u, v)$ is called a supersolution $($subsolution$)$ of system \eqref{colv} in $(0,+\i)\times\bar\Omega$.
Furthermore, if  both $({u}_1, {v}_1)$ and $({u}_2, {v}_2)$ are supersolutions $($subsolutions$)$ of \eqref{colv} in $(0,+\i)\times\bar\Omega$,
 then $\min(({u}_1, {v}_1), ({u}_2, {v}_2))$ $(\max(( u_1, v_1), (u_2, v_2)))$ is also called a supersolution  $($subsolution$)$ of \eqref{colv} in $(0,+\i)\times\bar\Omega$,
 where min and max are to be understood componentwise.
\end{defi}

The following comparison principle is derived directly from \cite{FT,PM}.
\begin{lem}\label{com}
If $(\u u,\u v)$ and $(\bar u,\bar v)$ are sub- and supersolutions of \eqref{colv} respectively,
and it holds  $ (\u u,\u v)(0,\cdot)\leq(\bar u,\bar v)(0,\cdot)$ in $\Omega$, then
$(\u u,\u v)\leq(\bar u,\bar v)$ in $[0,+\i)\times\bar\Omega$.
\end{lem}

Denote
$$
f_1(u,v)=r_1u(1-a_1-u+a_1v)\ \text{ and }\
f_2(u,v)=r_2(1-v)(a_2u-v).
$$
Since ${\p f_1(u,v)}/{\p v}\geq0$ and ${\p f_2(u,v)}/{\p u}\geq0$,
one gets that
\begin{align*}
\begin{cases}
  u_t\geq d_1\Delta u+u\int_0^1 \frac{\p f_1}{\p u}(\tau u,\tau v)d\tau&\text{ in }\ (0,+\i)\times\Omega, \\
v_t\geq d_2\Delta v+v\int_0^1 \frac{\p f_2}{\p v}(\tau u,\tau v)d\tau&\text{ in }\ (0,+\i)\times\Omega.
\end{cases}
\end{align*}
By \cite[Theorem 3.1]{G}, one immediately obtains the following weak Harnack inequalities.
\begin{lem}\label{harnack}
Given any $(t_0,x_0)\in(0,+\i)\times\Omega$, any $0\leq\theta_1<\theta_2<\theta\leq1$ and any $\alpha,\beta\in(0,1)$.
Let  $R\in(0,1]$ be such that $B(x_0,R)\subset\Omega$. There exist positive constants $p=p(N,d_1,\theta_1,\theta_2,\alpha,\beta)$,
$q=q(N,d_2,\theta_1,\theta_2,\alpha,\beta)$, $C=C(R,N,d_1,\theta_1,\theta_2,\alpha,\beta)$
and $\tilde C=\tilde  C(R,N,d_2,\theta_1,\theta_2,\alpha,\beta)$  such that
$$
\|u\|_{L^{p}((t_0,t_0+\theta_1 R^2)\times B(x_0,\alpha R))}
\leq C\inf_{(t,x)\in(t_0+\theta_2R^2,t_0+\theta R^2)\times B(x_0,\beta R)}u(t,x)
$$
and
$$
\|v\|_{L^{q}((t_0,t_0+\theta_1 R^2)\times B(x_0,\alpha R))}
\leq\tilde C\inf_{(t,x)\in(t_0+\theta_2R^2,t_0+\theta R^2)\times B(x_0,\beta R)}v(t,x).
$$
\end{lem}

When applying above conclusions in the following sections, we need to transform the cooperation-diffusion system \eqref{colv}
into the competition-diffusion system \eqref{clv}.
This transformation is achieved by changing variables as follows: $u=u$, $v=1-v$, $u_0=u_0$ and $v_0=1-v_0$.

\section{Fundamental properties}\label{S1}
This section is devoted to  fundamental properties of the spreading speeds.
In Section \ref{s12}, we study relationships between the spreading speeds and their dependence on position $z$, that is, we prove  formulas \eqref{dzz}
and \eqref{dzzz}
in Proposition \ref{dz} and Theorem \ref{idz}.
In Section \ref{s13}, we focus on the dependence of the spreading speeds on initial values, that is, we prove Theorem \ref{idi}.
In Section \ref{s21},  we study the upper bounds of the spreading speeds on a large class of  domains satisfying the extension property, that is, we prove Theorem \ref{upd}.

\subsection{Relationships between spreading speeds}\label{s12}

\begin{pr}[Proof of formulas \eqref{dzz} and \eqref{dzzz} in Proposition \ref{dz}]
  Let $\Omega\subset\mathbb R^N$ be any unbounded domain in a given direction $e\in\mathbb S^{N-1}$.
  Let $R=R(e)$ be the nonnegative real number  defined as in Definition \ref{unb}.
  Let $(u,v)$ be the solution of \eqref{colv} with any given initial value $(u_0,v_0)\in\Theta$.
  We only need to prove \eqref{dzz}
 since \eqref{dzzz} can be obtained by the same  arguments.
 From \eqref{w}, it is evident that  \eqref{dzz} holds when $w^*(e,u_0,v_0)=0$.
  Assume without loss of generality that $w^*(e,u_0,v_0)>0$.
 Fix any $\vp\in(0,w^*(e,u_0,v_0))$ and set
 $$
 \gamma=w^*(e,u_0,v_0)-\vp.
 $$
Owing to Definition \ref{s1}, there exists $A>R$ such that
$$
\sup_{s\geq\gamma t,\ x\in \bar{B(se,A)}\cap\bar\Omega}\left(|u(t,x)|+|1-v(t,x)|\right)\not\to0
\ \text{ as }\ t\to+\i.
$$
Hence, there exist some sequence $(t_n)_{n\in\mathbb N}\subset\R$ such that $t_n\to+\i$ as $n\to+\i$, $(s_n)_{n\in \mathbb N}\subset\R$
with $s_n\geq\gamma t_n$ and  $(x_n)_{n\in \mathbb N}\subset\bar{B_A}$ with
$x_n+s_n e\in\bar\Omega$ such that
\begin{align}\label{uvn}
\liminf_{n\to+\i}(|u(t_n,x_n+s_ne)|+|1-v(t_n,x_n+s_ne)|)>0,
\end{align}
where $B_A$ is the open Euclidean ball in $\R^N$ with the origin as its center and $A$ as its radius.
Up to extraction of a subsequence, we can assume that $x_n\to z\in \bar{B_A}$ as $n\to+\i$.

Assume on contrary that $w^*(e,z,u_0,v_0)<\gamma$.  By Definition \ref{s2}, one infers that
there exists $A'>0$ such that
$$
\sup_{s\geq\gamma t,\ x\in B_{A'},\ x+z+se\in\bar\Omega}(|u(t,x+z+se)|+|1-v(t,x+z+se)|)\to0\text{ as }t\to+\i.
$$
Note that $s_n\geq \gamma t_n$ and $(x_n-z)+z+s_ne=x_n+s_ne\in\bar\Omega$ for each $n\in\mathbb N$.
Since $x_n-z\in B_{A'}$ for all large $n$, one obtains that
$|u(t_n,x_n+s_ne)|+|1-v(t_n,x_n+s_ne))|\to0$ as $n\to+\i$. This is a contradiction with \eqref{uvn}.

As a conclusion, $w^*(e,z,u_0,v_0)\geq w^*(e,u_0,v_0)-\vp$ for all $\vp>0$. Together with \eqref{w} and arbitrariness of $\vp$, the  formula \eqref{dzz} follows. The proof is complete.
\end{pr}

\vspace{0.1cm}

\begin{pr}[Proof of Theorem \ref{idz}] We only need to prove \eqref{idzz},
the formula \eqref{idzzz} can be obtained similarly.
Assume that $\Omega$ is a strongly unbounded domain in a given direction $e$,
that there exists $\vp_0\in(0,1)$ such that $\Omega$ satisfies the interior ball property with a radius $\vp_0$,
  and that hypothesis $H_{y,z}$ is satisfied for some points $y$ and $z$ in $\R^N$. Then there exist
 $R_y>R(e,y)$, $R_z>R(e,z)$, $A>0$ and $s_0>0$ such that
\begin{align}\label{hyz1}
\forall\ s\geq s_0,\ \sup_{\substack{y'\in\bar{B(y+se,R_y)}\cap\bar\Omega\\
z'\in\bar{B(z+se,R_z)}\cap\bar\Omega}}d_\Omega(y',z')\leq A.
\end{align}
By Remark \ref{rks}, there holds
 $$
 \forall\ s\geq s_0,\ \bar{B(y+se,R_y)}\cap\bar\Omega\neq\emptyset\
 \text{ and }\
\bar{B(z+se,R_z)}\cap\bar\Omega\neq\emptyset.
 $$
Even if it means decreasing $\vp_0$, one can assume that $\vp_0<R_y$.

Let $(u,v)$ be the solution of \eqref{clv} in $[0,+\i)\times\bar\Omega$ with any given initial value $(u_0,v_0)\in\Theta$.
 In the case that both spreading speeds $w^*(e,y,u_0,v_0)$ and $w^*(e,z,u_0,v_0)$ are infinite,
 it is evident that $w^*(e,y,u_0,v_0)=w^*(e,z,u_0,v_0)$.
 Assume without loss of generality that $w^*(e,z,u_0,v_0)$ is finite.
Fix any $c>w^*(e,z,u_0,v_0)$.
By Definition \ref{s2}, there exist $c'\in(w^*(e,z,u_0,v_0),c)$ and $R>0$ such that
$$
\sup_{s'\geq c't',\ z'\in \bar{B(z+s'e, R)}\cap\bar\Omega} \left(|u(t',z')|+|1-v(t',z')|\right)\to0\
\text{ as }\ t'\to+\i.
$$
From Remark \ref{rks}, even if it means decreasing $R$, one can assume that $R\leq R_z$.
Let $\vp$ be any positive constant. Then there exists $t_0\geq1$ such that
$$
\forall\ t'\geq t_0,\ \forall\ s'\geq c't',\ \sup_{z'\in\bar{ B(z+s'e, R)}\cap\bar\Omega}\e(|u(t',z')|+|1-v(t',z')|\r)\leq\vp.
$$
By increasing $t_0$ if necessary, one can assume that $ct_0\geq s_0$ and $ct\geq c'(t+2)$ for all $t\geq t_0$.
Hence, there holds
\begin{align}\label{aa}
\forall\ t\geq t_0,\ \forall\ s\geq ct,\   \sup_{\tau\in(t,t+2),\ z'\in\bar {B(z+se, R)}\cap\bar\Omega}
\e(|u(\tau,z')|+|1-v(\tau,z')|\r)\leq\vp.
\end{align}

Fix any $t\geq t_0$ and any $s\geq ct$.
Define
 $$
 S_s=\{y'\in \bar{B(y+se,R_y)}\cap\bar\Omega:B(y',\vp_0)\subset{\Omega}\}.
 $$
Then
 \begin{align}\label{S}
  \bar{B(y+se,R_y)}\cap\bar\Omega\subset \bigcup_{y'\in S_s}\bar {B(y',\vp_0)}.
 \end{align}
Let $y'$ be any point in the set $S_s$ and let
$\beta\in(1/2,1)$ be any  constant.
 Since $\Omega$ satisfies the interior ball property with a radius $\vp_0$, with the aid of \eqref{hyz1},
 there exist $k$ points $y_1',\cdots,y_k'$ (see Figure \ref{fig3-4})
 \begin{figure}[htbp]
    \centering
    \begin{subfigure}[b]{0.49\textwidth}
        \centering
        \includegraphics[width=8cm]{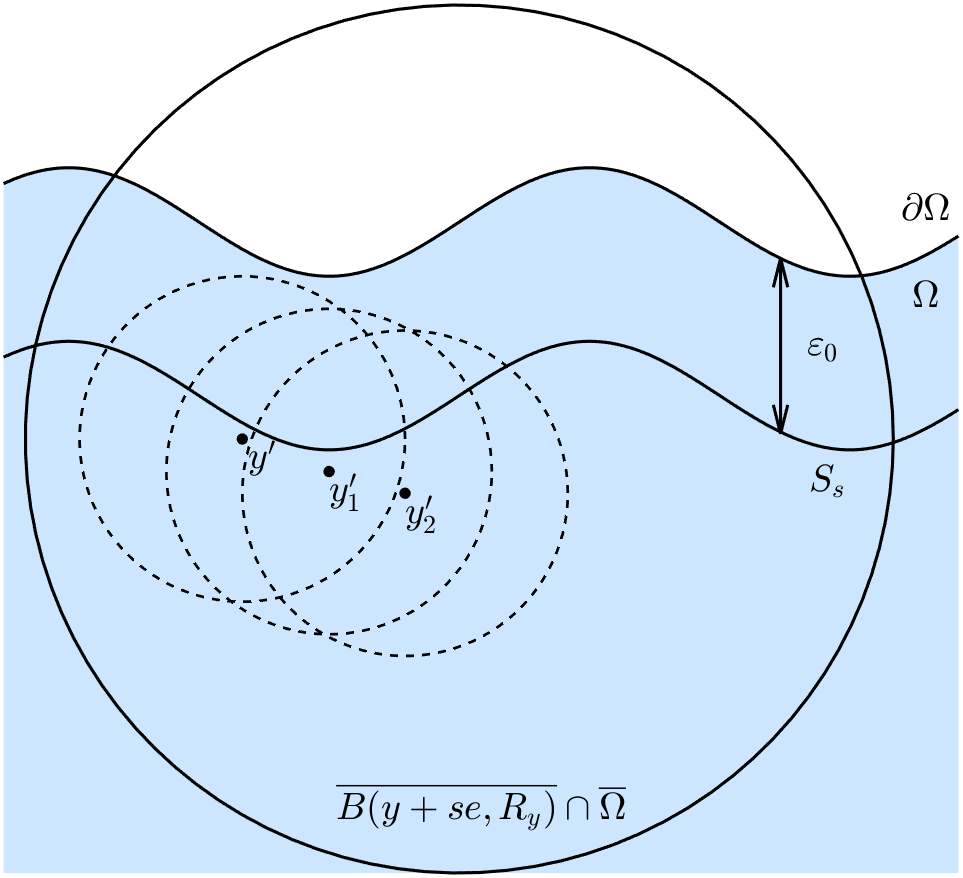}
    \end{subfigure}
    \hfill
    \begin{subfigure}[b]{0.49\textwidth}
        \centering
        \includegraphics[width=8cm]{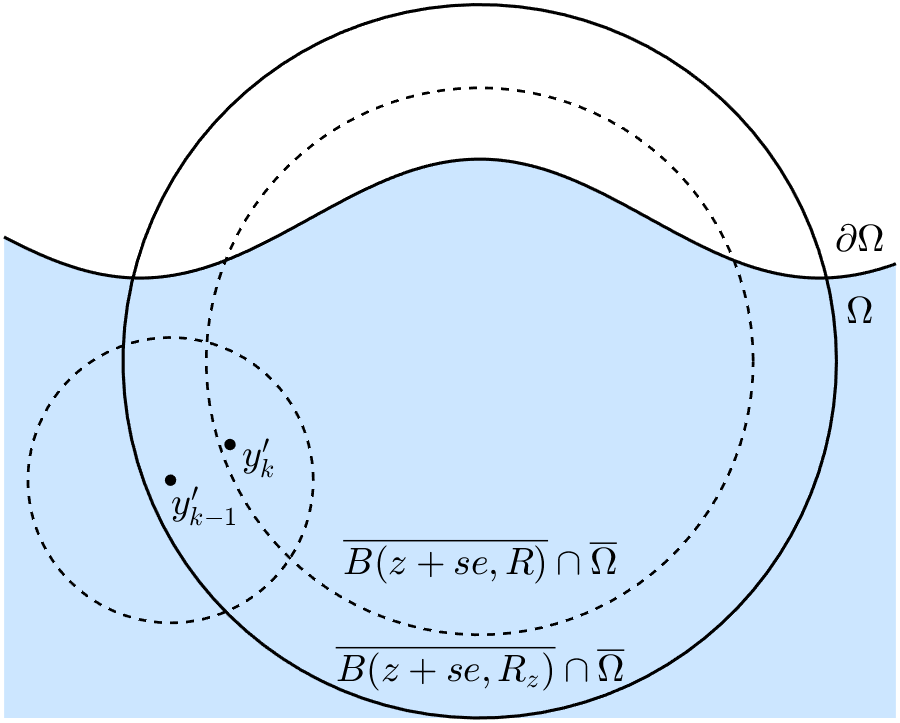}
    \end{subfigure}
    \caption{the choice of $k$ points $y_1',\cdots,y_k'$, where the shaded area is $\Omega$.}\label{fig3-4}
\end{figure}
such that
$$
\bc
y_1'\in B(y',\beta\vp_0),\ B(y_1',\vp_0/8)\subset B(y',\beta\vp_0)\not\supset \ B(y_1',\beta\vp_0),\\
B(y_i',\beta\vp_0)\subset\Omega, &1\leq i\leq k,\\
y_{i}'\in B(y_{i-1}',\beta\vp_0),\ B(y_i',\vp_0/8)\subset B(y_{i-1}',\beta\vp_0)\not\supset B(y_{i}',\beta\vp_0),
&
2\leq i\leq k,\\
y_k'\in \bar{B(z+se,R)}\cap\bar{\Omega}.
\ec
$$
By \eqref{hyz1}, the integer $k$ can be chosen to depend only on $A$.

For each $i=1,\cdots,k$, denote
$$
 t_i=t+\frac{i\vp_0^2}{4k}.
$$
 Clearly,  $t_k\in( t,t+1)$.
By Lemma \ref{harnack}, there exist two positive constants $p_0=p_0(N,d_1,\beta,k)$ and
$C_0=C_0(N,d_1,\beta,\vp_0,k)$ such that
\begin{align*}
\|u\|_{L^{p_0}\e(\e(t-\frac{\vp_0^2}{5k}, t+\frac{\vp_0^2}{5k}\r)\times B\e(y',\beta{\vp_0}\r)\r)}
\leq &C_0\inf_{(\tau,x)\in\e(t+\frac{\vp^2_0}{4k}, t+\frac{\vp_0^2}{2k}\r)\times B(y',\beta\vp_0)}u(\tau,x)\\
\leq& C_0\inf_{(\tau,x)\in\e(t_1, t_1+\frac{\vp_0^2}{5k}\r)\times B\e(y_1',\fr{\vp_0}{8}\r)}u(\tau,x).
\end{align*}
Similarly, for each $i=1,\cdots,k-1$,
there exist positive constants $p_i=p_i(N,d_1,\beta,k)$ and $C_i=C_i(N,d_1,\beta,\vp_0,k)$ such that
\begin{align*}
\|u\|_{L^{p_{i}}\e(\e(t_{i},t_{i}+\frac{\vp_{0}^2}{5k}\r)\times B\e(y_{i}',\fr{\vp_0}8\r)\r)}
&\leq C_{i}\inf_{(\tau,x)\in\e(t_{i+1}, t_{i+1}+\frac{\vp_{0}^2}{4k}\r)
\times B(y_{i}',\beta\vp_{0})}u(\tau,x)\\
&\leq C_{i}\inf_{(\tau,x)\in\e(t_{i+1}, t_{i+1}+\frac{ \vp_0^2}{5k}\r)
\times B\e(y_{i+1}',\fr{\vp_{0}}{8}\r)}u(\tau,x).
\end{align*}
In view of \eqref{aa} and the positivity of $u$ in $(0,+\i)\times\bar\Omega$ derived from \eqref{0uv1},
we conclude by induction that there is a positive constant $C=C(N,d_1,\beta,\vp_0,k)$ such that
\begin{align*}
\|u\|_{L^{p_0}\e(\e(t-\frac{\vp_0^2}{5k}, t+\frac{\vp_0^2}{5k}\r)\times B\e(y',\beta{\vp_0}\r)\r)}\notag
\leq\ & C_0
\left|\e(t_1, t_1+\frac{\vp_0^2}{5k}\r)\times B\e(y_1',\fr{\vp_0}{8}\r)\right|^{-\frac{1}{p_1}}
\|u\|_{L^{p_1}\e(\e(t_1, t_1+\frac{\vp_0^2}{5k}\r)\times B\e(y_1',\frac{\vp_0}{8}\r)\r)}\notag\\
\leq\ &C\inf_{(\tau,x)\in\e(t_k, t_{k}+\frac{ \vp_{0}^2}{5k}\r)
\times B\e(y_{k}',\frac{\vp_{0}}{8}\r)}u(\tau,x)\notag\\
\leq\ &C\sup_{(\tau,x)\in\e(t, t+2\r)\times (\bar{B(z+se,R)}\cap\bar\Omega)}u(\tau,x)\notag\\
\leq\ &C\vp .
\end{align*}
Similarly, there exist positive constants $q_0=q_0(N,d_2,\beta,k)$ and
$\tilde C=\tilde C(N,d_2,\beta,\vp_0,k)$  such that
 \begin{align*}
 \|1-v\|_{L^{q_0}\e(\e(t-\frac{\vp_0^2}{5k}, t+\frac{\vp_0^2}{5k}\r)\times B\e(y',\beta\vp_0\r)\r)}
\leq\tilde C\sup_{(\tau,x)\in\e(t, t+2\r)\times (\bar{B(z+se,R)}\cap\bar\Omega)}(1-v(\tau,x))\leq\tilde C \vp.
\end{align*}

For all $t\geq t_0$,
define the set
$$
Q_{t,\beta}=\left\{
(\tau,x)\in \e(t-\frac{\vp_0^2}{5k}, t+\frac{\vp_0^2}{5k}\r)\times B(y',\beta{\vp_0}):\
 y'\in S_s,\ s\geq ct\right\}.
$$
Since $\vp$  is arbitrary and  since the constants $p_0$, $q_0$, $C$ and $\tilde C$  are independent of $\vp$, one gets that
  \begin{align}\label{up0}
\|u\|_{L^{p_0}\e(Q_{t,\beta}\r)}+\|1-v\|_{L^{q_0}\e(Q_{t,\beta}\r)}\to0\ \text{ as }\ t\to+\i.
\end{align}
Assume on contrary that there exists $\gamma_0>0$  such that
  \begin{align*}
 \sup_{ (\tau,x)\in Q_{t,\beta}}
 |u(\tau,x)|\geq\gamma_0
 \ \text{ as }t\to+\i.
 \end{align*}
 For all large $t$, define the set
 $$
 D_{t,\beta}=\e\{(\tau,x)\in Q_{t,\beta}:\ |u(\tau,x)|\geq\frac{\gamma_0}{2}\r\}.
 $$
 Since $u_t$ and $\nabla u$ are globally bounded in $(0,+\i)\times\bar\Omega$ by standard parabolic estimates (see \cite{lib})
 and since $u\in(0,1)$ is continuous in $(0,+\i)\times\bar\Omega$,
one obtains that there exists a constant $\sigma_0>0$ such that
$|D_{t,\beta}|>\sigma_0$ for all large $t$.
  Hence, for all large $t$,
   \begin{align*}
\|u\|_{L^{p_0}\e(Q_{t,\beta}\r)}\geq\|u\|_{L^{p_0}(D_{t,\beta})} \geq\frac{\gamma_0\sigma_0^{1/p_0}}{2},
\end{align*}
it is a contradiction with \eqref{up0}.
Applying similar arguments to $v$, one has
 \begin{align}\label{pppp}
 \sup_{(\tau,x)\in Q_{t,\beta}}
 \e(|u(\tau,x)|+|1-v(\tau,x)|\r)\to0\ \text{ as }\ t\to+\i.
 \end{align}

Letting $\beta\to1$, the continuity of $(u,v)$ gives rise to
 \begin{align*}
\sup_{  y'\in S_s, \ s\geq ct, \ (\tau,x)\in\e(t-\frac{\vp_0^2}{5k},t+\frac{\vp_0^2}{5k}\r)\times B\e(y',\vp_0\r)
 }\e(|u(t,x)|+|1-v(t,x)|\r)\to0\ \text{ as }\ t\to+\i.
 \end{align*}
 For any $x\in\p B(y',\vp_0)$, there exists a sequence of points $(y_n')_{n\in\mathbb N}\subset B(y',\vp_0)$
 such that $y_n'\to x$ as $n\to+\i$.
Since $u,v\in C^{0,1}([0,+\i)\times\bar\Omega)$, one obtains that
 \begin{align*}
\sup_{   y'\in S_s,\ s\geq ct,\ (\tau,x)\in\e(t-\frac{\vp_0^2}{5k},t+\frac{\vp_0^2}{5k}\r)\times\bar {B(y',\vp_0)}
 }\e(|u(t,x)|+|1-v(t,x)|\r)\to0\ \text{ as }\ t\to+\i.
 \end{align*}
According to \eqref{S}, one has
 \begin{align*}
\sup_{  s\geq ct, \ y'\in \bar{B(y+se,R_y)}\cap\bar\Omega }\e(|u(t,y')|+|1-v(t,y')|\r)\to0\ \text{ as }\ t\to+\i.
 \end{align*}
 Owing to Definition \ref{s2}, one gets that $w^*(e,y,u_0,v_0)$ is finite and $w^*(e,y,u_0,v_0)\leq c$.
 Since this inequality holds for all $c>w^*(e,z,u_0,v_0)$, one obtains that $w^*(e,y,u_0,v_0)\leq w^*(e,z,u_0,v_0)$.

 By changing the roles of $y$ and $z$, one gets that $w^*(e,y,u_0,v_0)= w^*(e,z,u_0,v_0)$. The proof is complete.
\end{pr}

\subsection{Dependence on initial values}\label{s13}

\begin{pr}[Proof of Theorem \ref{idi}]
The proof is divided into two steps.

{\bf Step 1.} Let $(u,v)$ be the solution of \eqref{clv} with a given initial value $(u_0,v_0)\in \Theta$.
  We shall prove that $(u,v)(t,x)\to(1,0)$ locally uniformly in $x\in\bar\Omega$ as $t\to+\i$.
For all $r>0$ and $z\in\R^N$, define
\begin{align}\label{ld}
\lambda_r^z=\inf_{\psi\in C_c^1(\bar\Omega\cap B(z,r)),\ \psi\not\equiv0}
\frac{\int_{x\in\Omega\cap B(z,r)}|\nabla\psi(x)|^2dx}{\int_{x\in\Omega\cap B(z,r)}\psi^2(x)dx},
\end{align}
where $C_c^1(\bar\Omega\cap B(z,r))$ denotes the set of functions which are of class $C^1$ in $\bar\Omega\cap B(z,r)$
and are compact supported in $\bar\Omega\cap B(z,r)$.
 From \cite[Lemma 2.2]{bhn}, there exists $R\geq\max( R_0,1)$ such that
 \begin{align}\label{ldz}
 \forall\ z\in\R^N,\ \lambda_{2R}^z\leq \frac{r_1(1-a_1)}{2d_1}
 .
 \end{align}
 According to the assumption of Theorem \ref{idi}, one has
 $\Omega\cap B(z,r)\neq\emptyset$ for all $z\in\R^N$ and $r\geq R$.
 Let $u_{\rm KPP}$ be the solution of
 $$
 \bc
 (u_{\rm KPP})_t=d_1\Delta u_{\rm KPP}+r_1 u_{\rm KPP}(1-a_1-u_{\rm KPP})&\text{ in }\ (0,+\i)\times \Omega,\\
 (u_{\rm KPP})_\nu=0&\text{ on }\ (0,+\i)\times\p\Omega,\\
 u_{\rm KPP}(0,\cdot)=\tilde u_0 &\text{ in }\ \Omega,
 \ec
 $$
where $\tilde u_0\in C(\bar\Omega,[0,1])\backslash\{0\}$
is compact supported in $\bar\Omega$.
The maximum principle yields that  $0\leq u_{\rm KPP}\leq1$ in $[0,+\i)\times\bar\Omega$.

Recall that there exists $\vp_0\in(0,1)$ such that the connected domain $\Omega$
satisfies the interior ball property with a radius $\vp_0$.
For each $z\in\R^N$ and $r\geq R$, denote
$$
S_{z,r}=\{z'\in\bar\Omega\cap\bar{ B(z,r)}:B(z',\vp_0)\subset\Omega\}.
$$
Then $\bar\Omega\cap\bar{ B(z,r)}\subset\bigcup_{z'\in S_{z,r}}\bar{B(z',\vp_0)}$.
Choose $\alpha\in(1/2,1)$ is a constant sufficiently  close to $1$.
For any $z'\in S_{z,2R}$, it follows from Lemma \ref{harnack} that
 there exist positive constants $p=p(d_1,N,\alpha)$ and $C=C(d_1,N,\vp_0,\alpha)$ such that
 for all $t\geq1$, $ z\in\R^N$ and $z'\in S_{z,2R}$,
 \begin{align}\label{h3}
 C \| u_{\rm KPP}\|_{L^p\e(\e(t-\frac{ \vp_0^2}{5},t+\frac{ \vp_0^2}{5}\r)\times B(z',\alpha \vp_0)\r)}
\leq\inf_{(\tau,x)\in\e(t+\frac{ \vp_0^2}{4},t+\frac{ \vp_0^2}{2}\r)\times B(z', \alpha \vp_{0})} u_{\rm KPP}(\tau,x).
 \end{align}

Choose a constant $s_0>0$ sufficiently small.
Since $u(1,\cdot)$ is a positive continuous function in $\bar\Omega$, there exists
 $\tilde u_0\in C(\bar\Omega,[0,1])\backslash\{0\}$ such that
$$
{\rm supp} (\tilde u_0)\subset\bar\Omega\cap\bar{ B_R}\ \text{ and }\
\tilde u_0(x)\leq\min\left(1,Cs_0,u(1,x)\right)\text{ for }x\in\bar\Omega,
$$
where  $B_R$ is the open Euclidean ball in $\R^N$ with the origin as its center and $R$ as its radius,
${\rm supp}(\tilde u_0)$
denotes the support of $\tilde u_0$.

We now claim that
\begin{align}\label{c1}
\forall\ z_0\in\R^N,\ \exists\ t_{z_0}\geq1+\frac{\vp_0^2}{4},\ \  u_{\rm KPP}(t_{z_0},\cdot)\geq Cs_0
\ \text{ in }\ \bar\Omega\cap\bar{B(z_0,2R)}.
\end{align}
Assume that this is not true for some point $z_0\in\R^N$.
From \eqref{h3}, for all $z_0'\in S_{z_0,2R}$
and   $t\geq1+\vp_0^2/4$, there holds
 \begin{align*}
\|u_{\rm KPP}\|_{L^p\e(\e(t-\frac{\vp_0^2}{5},t+\frac{\vp_0^2}{5}\r)\times B(z_0', \alpha\vp_{0})\r)}
\leq C^{-1}\inf_{(\tau,x)\in\e(t+\frac{ \vp_0^2}{4},t+\frac{ \vp_0^2}{2}\r)\times( B(z_0', \alpha \vp_{0})\cap B(z_0,2R))} u_{\rm KPP}(\tau,x)
\leq s_0.
\end{align*}
 Since $(u_{\rm KPP})_t$ and  $\nabla u_{\rm KPP}$ are globally bounded in $[0,+\i)\times\bar\Omega$
 by standard parabolic estimates,
and since $ u_{\rm KPP}$ is continuous and positive in $(0,+\i)\times\bar\Omega$, even if it means decreasing $s_0$,
it follows from similar arguments as the proof of \eqref{pppp} that
\begin{align*}
\forall\ t\geq1+\frac{\vp_0^2}{4},\  \forall\ z_0'\in S_{z_0,2R},\
\sup_{(\tau,x)\in\e(t-\frac{\vp_0^2}{5},t+\frac{\vp_0^2}{5}\r)\times B(z_0',\alpha \vp_0)} u_{\rm KPP}(\tau,x)
\leq \frac{1}{4}(1-a_1).
 \end{align*}
By increasing $\alpha$ if necessary, one infers from the continuity of $u_{\rm KPP}$ and the definition  of  $S_{z_0,2R}$ that
 \begin{align}\label{tuv}
\forall\ t\geq1+\frac{\vp_0^2}{4},\
\sup_{x\in\bar{B(z_0,2 R)}\cap\bar\Omega}u_{\rm KPP}(t,x)
\leq \frac{1}{2}(1-a_1).
 \end{align}

By \eqref{tuv}, one knows that
\begin{align}\label{tuv1}
(u_{\rm KPP})_t(t,x)\geq d_1\Delta u_{\rm KPP}(t,x)+\frac{1}{2}r_1(1-a_1)u_{\rm KPP}(t,x)
\end{align}
for all $t\geq1+\vp_0^2/4$ and $x\in\bar{B(z_0, 2R)}\cap\bar\Omega$.
In term of \eqref{ld} and \eqref{ldz}, there is a function $w\in C_c^1(\bar\Omega\cap B(z_0, 2R))$
such that $w\not\equiv0$ and
$$
\int_{x\in\Omega\cap B(z_0,2R)}|\nabla w(x)|^2dx
<\frac{r_1(1-a_1)}{2d_1}\int_{x\in\Omega\cap B(z_0,2R)} w^2(x)dx.
$$
For all $t\geq1+\vp_0^2/4$, set
$$
\Lambda(t)=\int_{x\in\Omega\cap B(z_0,2R)}w^2(x)\ln \e(u_{\rm KPP}(t,x)\r)dx.
$$
Multiply \eqref{tuv1} by $w^2/u_{\rm KPP} \geq0$ and integrate by parts over $x\in\Omega\cap B(z_0,2R)$.
Since $\e(u_{\rm KPP}\r)_\nu=0$ on $\p\Omega$ and $w\in C_c^1(\bar\Omega\cap B(z_0, 2R))$, one gets that
\begin{align*}
\Lambda'(t)
\geq\ &\frac{1}{2}r_1(1-a_1)\int_{x\in\Omega\cap B(z_0,2R)}w^2(x)dx
-d_1\int_{x\in\Omega\cap B(z_0,2R)}\nabla u_{\rm KPP}(t,x)\cdot\nabla\left(\frac{w^2(x)}{u_{\rm KPP}(t,x)}\right)dx\\
=\ &\frac{1}{2}r_1(1-a_1)\int_{x\in\Omega\cap B(z_0,2R)}w^2(x)dx
   -2d_1\int_{x\in\Omega\cap B(z_0,2R)}\frac{w(x)\nabla u_{\rm KPP}(t,x)\cdot\nabla w(x)}{u_{\rm KPP}(t,x)}dx\\
   &+ d_1\int_{x\in\Omega\cap B(z_0,2R)}\frac{w^2(x)|\nabla u_{\rm KPP}(t,x)|^2}{u_{\rm KPP}^2(t,x)}dx\\
\geq \ &\frac{1}{2}r_1(1-a_1)\int_{x\in\Omega\cap B(z_0,2R)}w^2(x)dx
   -d_1\int_{x\in\Omega\cap B(z_0,2R)}|\nabla w(x)|^2dx\\
  >\ &0
\end{align*}
 for all $t\geq1+\vp_0^2/4$.
Then $\Lambda(t)\to+\i$ as $t\to+\i$. But it follows from \eqref{tuv} and $a_1\in(0,1)$ that
$$
\forall\ t\geq1+\frac{\vp_0^2}{4},\ \ \Lambda(t)\leq \ln\e(\frac{1-a_1}2\r)\times\int_{x\in\Omega\cap B(z_0,2R)} w^2(x)dx<0,
$$
it is  a contradiction. Therefore, the claim \eqref{c1} is true.

From the choice of $\tilde u_0$
and from \eqref{c1} applied at $z_0=0$ (the origin in $\R^N$), one gets that
$$
u_{\rm KPP}(t_0,\cdot)\geq \tilde u_0 \ \text{ in }\bar\Omega.
$$
The comparison principle implies that $u_{\rm KPP}(t+t_0,x) \geq u_{\rm KPP}(t,x)$
for all $t\geq0$ and $x\in\bar\Omega$.
Let $z_0$ be any point in $\R^N$.
By $\vp_0<R$, one gets that $B(z_0',\vp_0)\subset\bar\Omega\cap\bar{ B(z_0,2R)}$ for all $z_0'\in S_{z_0,R}$.
Choose a constant $\beta\in(1/2,1)$.
From \eqref{c1}, there holds
\begin{align}\label{cs0}
\forall\ k\in\mathbb N,\ \forall\ z_0'\in S_{z_0,R},\ \
\inf _{B(z_0',\beta\vp_0)}u_{\rm KPP}(t_{z_0}+kt_0,\cdot)\geq Cs_0.
\end{align}
By Lemma \ref{harnack} and \eqref{cs0},
 there exist constants $C_0'=C_0'(d_1,N,\vp_0,\beta)$, $\tilde C_0'=\tilde C_0'(d_1,N,\vp_0,\beta)$ and $p_0'=p_0'(d_1,N,\beta)$ such that
\begin{align*}
\inf_{(\tau,x)\in\e( t_{0,k}, t_{0,k}+\frac{2\vp_0^2}{5}\r)\times B(z_0',\beta \vp_0)}u_{\rm KPP}(\tau,x)
&\ \geq C_0'\| u_{\rm KPP}\|_{L^{p_0'}\left(\left( t_{z_0}+kt_0-\frac{\vp_0^2}{5},  t_{z_0}+kt_0+\frac{\vp_0^2}{5}\right)\times B( z_0',\beta\vp_0)\right)}\\
&\ \geq C_0'\|u_{\rm KPP}(t_{z_0}+k t_0,\cdot)\|_{L^{p_0'}\left( B( z_0',\beta\vp_0)\right)}
\geq \tilde C_0'C s_0
\end{align*}
for all $k\in\mathbb N$ and $z_0'\in S_{z_0,R}$,
where
$t_{0,k}=t_{z_0}+kt_0+{2\vp_0^2}/{5}$.
Let
$n_0=\min\left\{n\in\mathbb N:n\geq {5 t_0}/{\vp_0^2}\right\}$.
For each $k\in\mathbb N$ and each $i=1,\cdots,n_0$, define
$$
t_{i,k}=t_{0,k}+\frac{i\vp_0^2}{5}.
$$
By Lemma \ref{harnack}, for each $i=1,\cdots,n_0$, there  exist constants $C_i'=C_i'(d_1,N,\vp_0,\beta)$, $\tilde C_i'=\tilde C_i'(d_1,N,\vp_0,\beta)$ and $p_i'=p_i'(d_1,N,\beta)$ such that
\begin{align*}
\inf_{(\tau,x)\in\e(t_{i,k},t_{i,k}+\frac{2\vp_0^2}{5}\r)\times B(z_0',\beta \vp_0)} u_{\rm KPP}(\tau,x)
&\geq C_i'\| u_{\rm KPP}\|_{L^{p_i'}\left(\left(t_{i-1,k}, t_{i-1,k}+\frac{\vp_0^2}{10}\right)\times B( z_0',\beta\vp_0)\right)}
\geq \tilde C_i'C s_0
\end{align*}
for all $k\in\mathbb N$ and $z_0'\in S_{z_0,R}$.
According to the definition of $n_0$, one concludes by induction  that there exists $\tilde C=\tilde C(d_1,N,\vp_0,\beta)$ such that
\begin{align*}
\forall\ k\in\mathbb N,\
\forall\ z_0'\in S_{z,R},\ \
\inf_{(\tau,x)\in\e(t_{z_0}+kt_0+\frac{2\vp_0^2}{5},t_{z_0}+(k+1)t_0+\frac{4\vp_0^2}{5}\r)\times B(z_0',\beta \vp_0)}u_{\rm KPP}(\tau,x)
\geq \tilde CC s_0.
\end{align*}
Therefore,  one gets that
\begin{align*}
\forall\  z_0'\in S_{z,R},\ \forall\ t\geq t_{z_0}+\frac{2\vp_0^2}{5},\ \
\inf_{B(z_0',\beta\vp_0)}u_{\rm KPP}(t,\cdot)
\geq \tilde CCs_0.
\end{align*}

By increasing $\beta$ if necessary,
the definition of $S_{z,R}$ and the continuity of $ u_{\rm KPP}$ lead to
\begin{align*}
\forall \ t\geq t_{z_0}+1,\  \forall\ z_0\in\R^N,\ \inf_{\bar{B(z_0,R)}\cap\bar\Omega}u_{\rm KPP}(t,\cdot)
\geq\frac{1}{2} \tilde CCs_0.
\end{align*}
From the choice of $\tilde u_0$ and the comparison principle, one gets that
$u(t+1,x)\geq u_{\rm KPP}(t,x)$ for all $t\geq0$ and $x\in\bar\Omega$.
Then, for all $t\geq t_{z_0}+2$ and $z_0\in\R^N$,
\begin{align}\label{uv}
\min_{\bar{B(z_0,R)}\cap\bar\Omega}u(t,\cdot)
\geq\frac{1}{2} \tilde CCs_0>0,
\end{align}
where $C$, $\tilde C$ and  $s_0$ do not depend on the point $z_0$.
Pick now any sequence of positive real numbers $(t_n)_{n\in\mathbb N}$ such that
$t_n\to+\i$ as $n\to+\i$.
Remember that $(0,0)\leq (u,v)\leq(1, 1)$ in $[0,+\i)\times\bar\Omega$.
From standard parabolic estimates, up to extraction of a subsequence, the functions
$$
(u_n,v_n)(t,x)=(u,v)(t+t_n,x)
$$
converge locally uniformly in $(t,x)\in\R\times\bar\Omega$ to a classical nonnegative solution $(U,V)$ of
$$
\bc
U_t=d_1\Delta U+r_1U(1-U-a_1V) &\text{ in }\ \R\times\Omega,\\
V_t=d_2\Delta V+r_2V(1-V-a_2U)&\text{ in }\ \R\times\Omega,\\
(U_\nu,V_\nu)=(0,0)&\text{ on }\ \R\times\p\Omega.
\ec
$$
It is evident that
$(U,1-V)\leq(1,1)$ in $\R\times\bar\Omega$.
By \eqref{uv}, one gets that
$$
 \forall\ t\in\R,\ \forall\ x\in \bar\Omega,\ \ U(t,x)\geq\frac{1}{2} \tilde CCs_0>0.
$$
Then,
 $0<\tilde CCs_0/2\leq1$.

Call now $(\omega_1(t),\omega_2(t))$ the solution of
\begin{align}\label{bd}
\bc
\dot\omega_1(t)=r_1\omega_1(t)(1-a_1-\omega_1(t)+a_1\omega_2(t)):=g_1(\omega_1(t)), &t>0,\\
\dot\omega_2(t)=r_2(1-\omega_2(t))(a_2\omega_1(t)-\omega_2(t)):=g_2(\omega_2(t)), &t>0,\\
(\omega_1(0),\omega_2(0))=(\tilde CCs_0/2,0).
\ec
\end{align}
Then $(0,0)\leq (\omega_1(t),\omega_2(t))\leq(1,1)$ for all $t\geq0$.
Denote $h(\omega_1,\omega_2)=\omega_1^{-1}(1-\omega_2)^{-1}$.
By a direct calculation, there holds
 $$
\frac{\p(g_1h)}{\p(\omega_1)}+\frac{\p(g_2h)}{\p(\omega_2)}=-r_1(1-\omega_2)^{-1}-r_2\omega_1^{-1}<0
$$
 for all $\omega_1\neq0$ and $\omega_2\neq1$.
 Owing to Bendixson-Dulac theorem, the system \eqref{bd} has no closed orbit.
Since $(1,1)$ is the unique stable equilibrium of \eqref{bd}, one gets that
\begin{align}\label{ode}
(\omega_1(t),\omega_2(t))\to(1,1)\ \text{ as }t\to+\i.
\end{align}

For any $t\in\R$ and any $T\geq0$, since $(U(t-T,\cdot),1-V(t-T,\cdot))\geq( \tilde CCs_0/2,0)$ in  $\bar\Omega$,
the comparison principle implies that $(U(t,\cdot),1-V(t,\cdot))\geq(\omega_1(T),\omega_2(T))$
 in $\bar\Omega$.
Since this holds for all $t\in\R$ and
$T\geq0$, one concludes that
$(U(t,x),1-V(t,x))\geq(1,1)$ for all $(t,x)\in\R\times \bar\Omega$.
As a conclusion, there holds
$$
(U,V)=(1,0)\ \text{ in }\R\times\bar\Omega.
$$
By uniqueness of the limit, it follows that
\begin{align}\label{wd}
(u,v)(t,x)\to(1,0)\  \text{ locally uniformly in }x\in\bar\Omega \text{ as }t\to+\i.
\end{align}

{\bf Step 2.} We now turn to prove that the spreading speeds are independent of the choice of initial values.
Let now $(u_0,v_0)$ and $(\tilde u_0,\tilde v_0)$ be two initial values in $\Theta$.
Assume that $u_0$, $\tilde u_0$, $1-v_0$ and $1-\tilde v_0$ are less than $1$.
Let $(\tilde u,\tilde v)$ be the solution of \eqref{clv} in $[0,+\i)\times\bar\Omega$ with initial value $(\tilde u_0,\tilde v_0)$.
Let $e$ be a unit vector in $\R^N$. Notice that the assumptions in Theorem \ref{idi} yield that
$\Omega$ is strongly unbounded in the direction $e$.
Since $\max_{\bar\Omega} \tilde u_0<1$ and  $\max_{\bar\Omega}(1- \tilde v_0)<1$,
 and since $(u_0,1-v_0)$ is compactly supported,
it follows from \eqref{wd} that
 there exists $t_0>0$ such that
 $$
 (u(t_0,x),1-v(t_0,x))\geq (\tilde u_0(x),1-\tilde v_0(x))
 $$
 for all $x\in\bar\Omega$.
 The comparison principle implies that $ (u(t+t_0,x),1-v(t+t_0,x))\geq (\tilde u(t,x),1-\tilde v(t,x))$
 for all $t\geq0$ and $x\in\bar\Omega$, whence
 $w^*(e,u_0,v_0)\geq w^*(e,\tilde u_0,\tilde v_0)$.

 Changing the roles of $(u,v)$ and $(\tilde u,\tilde v)$, one obtains that  $w^*(e,u_0,v_0)\leq w^*(e,\tilde u_0,\tilde v_0)$.
 As a conclusion,  $w^*(e,u_0,v_0)=w^*(e,\tilde u_0,\tilde v_0)$.
Notice that $w_*(e,u_0,v_0)$ and  $w_*(e,z,u_0,v_0)$ are well-defined by \eqref{wd}.
For all $e\in \mathbb S^{N-1}$, $z\in\R^N$ and $(u_0,v_0)\in\Theta$,
the same arguments  imply that $w_*(e,u_0,v_0)$, $w^*(e,z,u_0,v_0)$ and $w_*(e,z,u_0,v_0)$
also do not depend on
the initial value $(u_0,v_0)$ provided that $(u_0,1-v_0)\ll(1,1)$.
The proof of Theorem \ref{idi} is thereby complete.
\end{pr}

\subsection{Upper bounds for domains with the extension property}\label{s21}

\begin{pr}[Proof of Theorem \ref{upd}]
Let $\Omega$ be a locally $C^2$ connected open subset of $\R^N$ satisfying the extension property.
Assume that $\Omega$ is strongly unbounded in a direction $e\in\mathbb S^{N-1}$.
Let $(u,v)$ be the solution of \eqref{clv} in $[0,+\i)\times\bar\Omega$ with
 initial value  $(u_0,v_0)\in\Theta$.
 Let $R_0>0$ be such that $B_{R_0}$ contains the supports of $u_0$ and $1-v_0$,
where $B_{R_0}$ is the open Euclidean ball in $\R^N$ with the origin as its center and $R_0$ as its radius.

{\bf Step 1: proof of part (a).}
Assume that $(u,v)(t,x)\to(1,0)$ locally uniformly in $x\in\bar\Omega$ as $t\to+\i$.
It follows that the speeds $w_*(e,u_0,v_0)$ and $w_*(e,z,u_0,v_0)$ are well-defined.
Assume that if \eqref{updu} holds, then  there exists $z_0\in\R^N$ such that
$w_*(e,z_0,u_0,v_0)>2\sqrt{d_1r_1}$.
Choose positive constants $c'$ and $\vp'$  such that
$2\sqrt{d_1r_1}<c'-2\vp'<c'+2\vp'<w_*(e,z_0,u_0,v_0)$.
By Definition \ref{ss2}, there exists $A>0$ such that
$$
\limsup_{t\to+\i}
  \left(\sup_{ (c'-\vp')t\leq s\leq (c'+\vp')t,\ x\in\bar{B(z_0+se,A)}\cap\bar\Omega}
  \left(|1-u(t,x)|+|v(t,x)|\right)\right)=0.
$$
On the other hand, by \eqref{updu} and \eqref{0uv1}, one has
\begin{align*}
0\leq\limsup_{t\to+\i}
  \left(\sup_{\substack{(c'-\vp')t\leq s\leq (c'+\vp')t,\\ x\in\bar{B(z_0+se,A)}\cap\bar\Omega}}
  u(t,x)\right)
 & \leq
\limsup_{t\to+\i}
  \left(\sup_{\substack{(c'-2\vp')t\leq|x|\leq(c'+2\vp')t}}
  u(t,x)\right)\\
 & \leq
  \limsup_{t\to+\i}
  \left(\sup_{\substack{|x|\geq(c'-2\vp')t}}
  u(t,x)\right)=0,
\end{align*}
this is a contradiction.
As a conclusion,  \eqref{updd} is true if \eqref{updu} holds.

We only need to prove \eqref{updu}.
Let $u_{\rm KPP}$ be the solution of
$$
\bc
(u_{\rm KPP})_t=d_1 \Delta u_{\rm KPP}+r_1u_{\rm KPP}(1-u_{\rm KPP})
& \text{ in }\ (0,+\i)\times\Omega,\\
(u_{\rm KPP})_\nu=0
& \text{ on }\ (0,+\i)\times\p\Omega,\\
u_{\rm KPP}(0,\cdot)=u_0
& \text{ in }\ \Omega.
\ec
$$
By similar arguments as in the proof of \cite[Theorem 1.8]{bhn}, there holds
$$
\forall\ c>2\sqrt{d_1r_1},\ \
\sup_{|x|\geq ct,\ x\in\bar\Omega}u_{\rm KPP}(t,x)\to0\ \text{ as }\ t\to+\i.
$$
Furthermore,
the comparison principle implies that
$u\leq u_{\rm KPP}$ in $[0,+\i)\times\bar\Omega$.
Together with \eqref{0uv1}, the formula \eqref{updu} follows.

{\bf Step 2: proof of part (b).}
Notice that the property \eqref{upd2} is stronger than \eqref{upd1},
we just need to prove \eqref{upd2}.
Assume that $d_1=d_2$ and $a_1a_2\leq1$.
Fix a speed $c>2\sqrt{d_1r_1(1-a_1)}$.
Take $C_0>4$, $\vp_0>0$
and $t_0>0$ be such that
\begin{align}\label{heatk}
\forall\ t\geq t_0,\ \forall\ z\in B_{R_0},
\ \forall\ |x|\geq ct,\ \
\frac{|z-x|^2}{C_0d_1t}\geq (r_1(1-a_1)+\vp_0)t.
\end{align}

Call $w(t,x)$  the solution of
\begin{align*}
\bc
w_t=d_1\Delta w &\text{ in }\ (0,+\i)\times\Omega,\\
w_\nu=0  &\text{ on }\ (0,+\i)\times\p\Omega,\\
w(0,\cdot)=\max(u_0,1-v_0) &\text{ in }\ \Omega.
\ec
\end{align*}
The maximum principle implies that $0\leq w\leq1$ in $[0,+\i)\times\bar\Omega$.
 For all $(t,x)\in[0,+\i)\times\bar\Omega$, define the function
$$
(\bar u,\bar v)(t,x)=\min((w(t,x)e^{r_1(1-a_1)t},a_2w(t,x)e^{r_1(1-a_1)t}),(1,1)).
$$
We shall prove that $(\bar u,\bar v)$ is a supersolution of \eqref{colv} in $[0,+\i)\times\bar\Omega$.
By $a_2>1$, one gets that
$(\bar u,\bar v)(0,\cdot)\geq (u_0,1-v_0)$ in $\bar\Omega$.
It is evident that
$\bar u_\nu(t,x)=0$ for all $t>0$ and $x\in\p\Omega$
such that $\bar u(t,x)<1$,
and $\bar v_\nu(t,x)=0$ for all $t>0$ and $x\in\p\Omega$
such that $\bar v(t,x)<1$.
It suffices to check that
$
  \mathcal L_1(\bar u,\bar v)\geq0
$
and
$
  \mathcal L_2(\bar u,\bar v)\geq0
$
(the operators $\mathcal L_1$ and $\mathcal L_2$ are defined as in Definition \ref{dss})
in two cases: $(\bar u,\bar v)=(we^{r_1(1-a_1)t},a_2we^{r_1(1-a_1)t})$ and
$(\bar u,\bar v)=(we^{r_1(1-a_1)t},1)$.
In the first case,
 thanks to $a_1a_2\leq1$ and $a_2\bar u=\bar v$, one gets that
\begin{align*}
  \mathcal L_1(\bar u,\bar v)
  =-r_1(a_1a_2-1)w^2e^{2r_1(1-a_1)t}\geq0
\end{align*}
and
\begin{align*}
  \mathcal L_2(\bar u,\bar v)=a_2r_1(1-a_1)we^{r_1(1-a_1)t}\geq0.
\end{align*}
In the second case, there holds $a_2we^{r_1(1-a_1)t}\geq 1$. Clearly, $\mathcal L_2(\bar u,\bar v)=0$.
By   $a_1a_2\leq1$, one gets that
\begin{align*}
  \mathcal L_1(\bar u,\bar v)
  =r_1we^{r_1(1-a_1)t}(we^{r_1(1-a_1)t}-a_1)
  \geq r_1we^{r_1(1-a_1)t}(a_2^{-1}-a_1)\geq 0.
\end{align*}

The comparison principle yields that
$
(0,0)\leq (u,1- v)\leq(\bar u, \bar v)
$
for all $t\geq0$ and $x\in\bar\Omega$.
The function $w$ can be written as
$$
w(t,x)=\int _\Omega p(d_1t,z,x)w(0,z)dz
\leq \int _{B_{R_0}}p(d_1t,z,x)w(0,z)dz,
$$
where $p(t,z,x)$ is the heat kernel in $\bar\Omega$ with Neumann boundary condition.
Because $\Omega$ is a locally $C^2$ connected open subset in $\R^N$ satisfying the extension property,
from \cite[Propostion 2.5]{bhn},
there are two positive constants $C$ and $\dl$ such that
$$
w(t,x)\leq C\|w(0,\cdot)\|_{L^\i(\bar\Omega)}(1+(d_1\dl t)^{-N/2}) \int _{B_{R_0}}e^{ -\frac{d_\Omega(z,x)^2}{C_0d_1t}}dz
$$
for all $t>0$ and $x\in\bar\Omega$.
 One concludes from \eqref{heatk} that
\begin{align*}
(0,0)\leq (u(t,x),1-v(t,x))&\leq
 Ce^{r_1(1-a_1)t}\|w(0,\cdot)\|_{L^\i(\bar\Omega)}(1+(d_1\dl t)^{-N/2})\int _{B_{R_0}}e^{ -\frac{d_\Omega(z,x)^2}{C_0d_1t}}dz
\times (1,a_2)\\
&\leq C e^{-\vp_0 t} \|w(0,\cdot)\|_{L^\i(\bar\Omega)}(1+(d_1\dl t)^{-N/2})|B_{R_0}|
\times (1,a_2)
\end{align*}
for all $t\geq t_0$ and $x\in\bar\Omega$ such that $|x|\geq ct$,
where $|B_{R_0}|$ denotes the Lebesgue measure of the ball $B_{R_0}$.
By letting $t\to+\i$, the estimate \eqref{upd2} follows.
The proof is thereby complete.
\end{pr}

\section{Exterior domains and domains containing large half-cylinders}\label{S2}
We focus on estimates of the spreading speeds in exterior domains and domains containing large half-cylinders.
In Section \ref{s22}, we establish general upper and lower bounds  for the spreading speeds in exterior domains
by proving Theorems \ref{upd-ex} and \ref{ext}.
In Section \ref{S23},
we provide a general lower bound   for the upper spreading speeds in domains containing large half-cylinders,
that is, we prove  Proposition \ref{half-cyl}.

\subsection{Exterior domains}\label{s22}
\begin{pr}[Proof of Theorem \ref{upd-ex}]
Let $\Omega$ be an exterior domain in $\R^N$ of class $C^2$ and
let $(u,v)$ be the solution of \eqref{clv} with initial value $(u_0,v_0)\in \Theta$.
Let $e\in\mathbb S^{N-1}$ be any fixed direction.
It is clearly that $\Omega$ is strongly unbounded in direction $e$.
Together with Theorem \ref{idi}, one concludes that the spreading speeds
$w_*(e,u_0,v_0)$, $w_*(e,z,u_0,v_0)$, $w^*(e,z,u_0,v_0)$ and $w^*(e,u_0,v_0)$
are well-defined.
Since $\R^N\backslash\Omega$ and the supports of $u_0$ and $1-v_0$ are compact, there exists a constant $R_0>0$
such that
\begin{align}\label{R0}
{\rm supp}(u_0),\ {\rm supp}(1-v_0),\ \p\Omega\subset B_{R_0}.
\end{align}
We consider two cases $2d_1\geq d_2$ and $2d_1<d_2$, respectively.

{\bf Case 1: $ 2d_1\geq d_2$.}
In this case,  $2\sqrt{r_1}\max(\sqrt{d_1},\sqrt{d_2/2})=2\sqrt{d_1 r_1}$.
By virtue of \eqref{w}, we just need to prove $w^*(e,u_0,v_0)\leq 2\sqrt{d_1r_1}$.
For all $t\geq0$ and $x\in\bar\Omega$, define the function
$$
(\bar u,\bar v)(t,x)=\min\left(\e(e^{-\sqrt{\frac{r_1}{d_1}}(x\cdot e-R_0)+2r_1t},a_2e^{-\sqrt{\frac{r_1}{d_1}}(x\cdot e-R_0)+2r_1t}\r),(1,1)\right).
$$
We shall prove that $(\bar u,\bar v)$ is a supersolution of \eqref{colv} in $[0,+\i)\times\bar\Omega$.
Let us first verify the initial and boundary conditions.
If $x\in\bar\Omega\cap\bar{ B_{R_0}}$, then
$(\bar u,\bar v)(0,\cdot)=(1,1)\geq (u_0,1-v_0)$
 since $|x\cdot e|\leq R_0$ and $a_2>1$.
If $x\in\bar\Omega\backslash\bar{ B_{R_0}}$, it then follows from \eqref{R0} that
$(\bar u,\bar v)(0,\cdot)\geq(0,0)=(u_0,1-v_0)$.
 Thus, there holds
 $$
 (\bar u,\bar v)(0,\cdot)\geq(u_0,1-v_0)\ \text{ in }\bar\Omega.
 $$
 Furthermore, one infers from \eqref{R0} that
 $$
 (\bar u,\bar v)=(1,1)\geq (u,1-v)\ \text{ on }(0,+\i)\times\p\Omega.
 $$
 By $a_2>1$, it suffices to prove that
 $\mathcal L_1(\bar u,\bar v)\geq0$
 and
 $\mathcal L_2(\bar u,\bar v)\geq0$
(the operators $\mathcal L_1$ and $\mathcal L_2$ are defined as in Definition \ref{dss})
in two cases:
$$
(\bar u,\bar v)=\e(e^{-\sqrt{\frac{r_1}{d_1}}(x\cdot e-R_0)+2r_1t},a_2e^{-\sqrt{\frac{r_1}{d_1}}(x\cdot e-R_0)+2r_1t}\r)
\ \text{ and }\
(\bar u,\bar v)=\e(e^{-\sqrt{\frac{r_1}{d_1}}(x\cdot e-R_0)+2r_1t},1\r).
$$
In the first case, thanks to $0<\bar u,\bar v<1$, $\bar v=a_2\bar u$ and $ 2d_1\geq d_2$, one gets that
 $$
\mathcal L_1(\bar u,\bar v)
= 2r_1\bar u-r_1\bar u-r_1\bar u(1-a_1-\bar u+a_1\bar v)
 \geq0
 $$
 and
 $$
 \mathcal L_2(\bar u,\bar v)
 = 2r_1\bar v-\frac{d_2r_1}{d_1}\bar v-r_2(1-\bar v)(a_2\bar u-\bar v)
=r_1\bar v\left(2-\frac{d_2}{d_1}\right)
 \geq0.
 $$
In the second case, similar argument leads to
$\mathcal L_1(\bar u,\bar v)\geq0$
 and
 $\mathcal L_2(\bar u,\bar v)=0$.

By the comparison principle, one has
$(0,0)\leq(u,1-v)\leq(\bar u,\bar v)$
in $[0,+\i)\times\bar\Omega$.
Let $c> 2\sqrt{d_1r_1}$ be any fixed speed.
For all $A>0$, $t\geq0$, $s\geq ct$ and  $x\in\bar{B(se,A)}\cap\bar\Omega$, since $x\cdot e\geq s- A\geq ct-A$,
there holds
\begin{align*}
 (0,0)\leq(u,1-v)
  \leq(1,a_2)\times e^{-\sqrt{\frac{r_1}{d_1}}((c-2\sqrt{d_1r_1})t-A-R_0)}.
\end{align*}
As a conclusion,
$$
\forall\ A>0,\ \forall \ c> 2\sqrt{d_1r_1},\ \
\sup_{s\geq ct,\ x\in\bar{B(se,A)}\cap\bar\Omega}(|u(t,x)|+|1-v(t,x)|)\to0\ \text{ as }\ t\to+\i,
$$
which implies that $w^*(e,u_0,v_0)\leq 2\sqrt{d_1r_1}$ and \eqref{ue1} follows.

{\bf Case 2: $ 2d_1< d_2$.}
In this case,  $2\sqrt{r_1}\max(\sqrt{d_1},\sqrt{d_2/2})=\sqrt{2d_2 r_1}$.
For all $t\geq0$ and $x\in\bar\Omega$, define the function
$$
(\bar u,\bar v)(t,x)=\min\left(\e(e^{-\sqrt{\frac{2r_1}{d_2}}(x\cdot e-R_0)+2r_1t},a_2e^{-\sqrt{\frac{2r_1}{d_2}}(x\cdot e-R_0)+2r_1t}\r),(1,1)\right).
$$
Let us show that $(\bar u,\bar v)$ is a supersolution of \eqref{colv} in $[0,+\i)\times\bar\Omega$.
By similar arguments as in Case 1, one gets that
$(\bar u,\bar v)(0,\cdot)\geq(u_0,1-v_0)$ in $\bar\Omega$ and
 $
 (\bar u,\bar v)\geq (u,1-v)$ on $(0,+\i)\times\p\Omega$.
 It suffices to prove that
 $\mathcal L_1(\bar u,\bar v)\geq0$
 and
 $\mathcal L_2(\bar u,\bar v)\geq0$
in two cases:
$$
(\bar u,\bar v)=\e(e^{-\sqrt{\frac{2r_1}{d_2}}(x\cdot e-R_0)+2r_1t},a_2e^{-\sqrt{\frac{2r_1}{d_2}}(x\cdot e-R_0)+2r_1t}\r)
\ \text{ and }\
(\bar u,\bar v)=\e(e^{-\sqrt{\frac{2r_1}{d_2}}(x\cdot e-R_0)+2r_1t},1\r).
$$
In the first case, since $0\leq\bar u,\bar v\leq1$, $\bar v=a_2\bar u$ and $ 2d_1<d_2$, there hold
 $$
 \mathcal L_1(\bar u,\bar v)
= 2r_1\bar u-\frac{2d_1r_1}{d_2}\bar u-r_1\bar u(1-a_1-\bar u+a_1\bar v)
\geq r_1\bar u\e(1-\frac{2d_1}{d_2}\r)
 \geq0
 $$
 and
 $$
 \mathcal L_2(\bar u,\bar v)
 = 2r_1\bar v-2r_1\bar v-r_2(1-\bar v)(a_2\bar u-\bar v)
=0.
 $$
 In the second case, similar argument yields that
 $\mathcal L_1(\bar u,\bar v)\geq0$
 and
 $\mathcal L_2(\bar u,\bar v)=0$.

The comparison principle leads to
$(0,0)\leq(u,1-v)\leq(\bar u,\bar v)$
in $[0,+\i)\times\bar\Omega$.
By the same argument as in Case 1, one has
$$
\forall\ A>0,\ \forall \ c> \sqrt{2d_2r_1},\ \
\sup_{s\geq ct,\ x\in\bar{B(se,A)}\cap\bar\Omega}(|u(t,x)|+|1-v(t,x)|)\to0\ \text{ as }\ t\to+\i.
$$
Then $w^*(e,u_0,v_0)\leq \sqrt{2d_2r_1}$.
Together with \eqref{w}, \eqref{ue1} follows. The proof of Theorem \ref{upd-ex} is complete.
\end{pr}

\begin{pr}[Proof of Theorem \ref{ext}]
Let $\Omega$ be an exterior domain in $\R^N$ of class $C^2$ and
let $(u,v)$ be the solution of \eqref{clv} with initial value $(u_0,v_0)\in \Theta$.
It suffices to prove formula \eqref{ext2}.
By Theorem \ref{idi}, formula \eqref{ext2} holds if $c=0$.
Fix a speed $c\in(0,2\sqrt{d_1r_1(1-a_1)})$.
Let $u_{\rm KPP}$ solve
\begin{align*}
\bc
 (u_{\rm KPP})_t=d_1\Delta u_{\rm KPP}+r_1u_{\rm KPP}(1-a_1-u_{\rm KPP}) &\ \text{ in }\ (0,+\i)\times\Omega,\\
(u_{\rm KPP})_\nu=0, &\ \text{ on }\ (0,+\i)\times\p\Omega,\\
u_{\rm KPP}(0,\cdot)=\min(u_0,1-a_1) &\ \text{ in }\ \Omega.
\ec
\end{align*}
The maximum principle yields that $0\leq u_{\rm KPP}\leq 1-a_1$ in $[0,+\i)\times\bar\Omega$.
From \cite[Theorem 1.9]{bhn}, one has
\begin{align*}
\forall\ 0\leq c'<2\sqrt{d_1r_1(1-a_1)},\ \
\sup_{|x|\leq c't,\ x\in\bar\Omega}|1-a_1- u_{\rm KPP}(t,x)|\to 0\ \text{ as }\ t\to+\i.
\end{align*}
By the comparison principle, there holds $u\geq u_{\rm KPP}$ in $[0,+\i)\times\bar\Omega$. Hence,
\begin{align}\label{bu}
\forall\ 0\leq c'<2\sqrt{d_1r_1(1-a_1)},\ \
\liminf_{t\to+\i}\left(\inf_{|x|\leq c't,\ x\in\bar\Omega} u(t,x)\right)
\geq \frac{1-a_1}{2}.
\end{align}

Assume on contrary that \eqref{ext2} is not true.
Together with Theorem \ref{idi},
there exist sequences $(t_n)_{n\in\mathbb N}\subset\R$ and $(x_n)_{n\in\mathbb N}\subset\bar\Omega$
such that $|x_n|\to+\i$, $t_n\to+\i$ and
$|x_n|/t_n\to\hat c\in[0, c]$ as $n\to+\i$, and such that
\begin{align}\label{contra}
\limsup_{n\to+\i}(|1-u(t_n,x_n)|+|v(t_n,x_n)|)>0.
\end{align}
For each $n\in\mathbb N$, define
$$
(u_n,v_n)(t,x)=(u,v)(t+t_n,x+x_n)
$$
for all $t\in\R$ such that $t\leq t_n$ and $x\in\bar\Omega-x_n$.
By standard parabolic estimates, up to extraction of a subsequence,
there is a vector-value function $(u_\i,v_\i)$ such that
 $(u_n,v_n)\to(u_\i,v_\i)$ locally uniformly in $\R\times\R^N$ as $n\to+\i$, where $(u_\i,v_\i)$ is the
  classical solution of
  \begin{align*}
\begin{cases}
  (u_\i)_t =d_1\Delta u_\i +r_1u_\i (1-u_\i -a_1v_\i )&\text{ in }\ \R\times\R^N, \\
(v_\i)_t =d_2\Delta v_\i +r_2v_\i (1-v_\i -a_2u_\i )&\text{ in } \ \R\times\R^N.
\end{cases}
\end{align*}
Thanks to \eqref{bu}, there holds
$$
\forall\ t\in\R, \ \forall\ x\in\R^N ,\ \
u_\i(t,x)\geq\frac{1-a_1}{4}.
$$

For any $t_0\in\R$, let $(\omega_1(t),\omega_2(t))$ be the solution of
$$
\bc
\dot \omega_1(t)=r_1\omega_1(t)(1-a_1-\omega_1(t)+a_1\omega_2(t)), &t>t_0,\\
\dot \omega_2(t)=r_2(1-\omega_2(t))(a_2\omega_1(t)-\omega_2(t)), &t> t_0,\\
(\omega_1(t_0),\omega_2(t_0))=\e(\frac{1-a_1}{4},0\r).
\ec
$$
Using the comparison principle, one has
$
(1,1)\geq (u_\i(t,x),1-v_\i(t,x))\geq (\omega_1(t+t_0),\omega_2(t+t_0))$
 for all $t\geq0$ and $x\in\R^N$.
Letting $t_0\to+\i$, there holds $(u_\i,1-v_\i)=(1,1)$ in $[0,+\i)\times\R^N$ since $(\omega_1(t),\omega_2(t))\to(1,1)$ as $t\to+\i$
(see the proof of \eqref{ode}).
Then one reaches a contradiction with \eqref{contra}.

Thus, we conclude that formula \eqref{ext2} holds.
The proof of Theorem \ref{ext} is thereby complete.
\end{pr}

\subsection{Domains containing large half-cylinders}\label{S23}

\begin{pr}[Proof of Proposition \ref{half-cyl}]
For any $R>0$, let $(\lambda_R,\psi_R)$ be the pair of first eigenvalue and first eigenfunction of
\begin{align*}
\bc
-d_1\Delta\psi_R=\lambda_R\psi_R &\text{ in }\ B_R,\\
\psi_R>0 &\text{ in }\ B_R,\\
\psi_R=0&\text{ on }\ \p B_R,\\
\|\psi_R\|_{L^\i(B_R)}=1,
\ec
\end{align*}
where $B_R$ is the open Euclidean ball in $\R^N$ with the origin as its center and $R$ as its radius.
Fix $\vp\in(0,2\sqrt{d_1r_1(1-a_1)}]$.
Since
$\lambda_R\to0$ as $R\to+\i$ (see \cite[Lemma 2.2]{bhn}), there exists $R_0>0$ large enough so that
\begin{align*}
\forall\ R\geq R_0,\ \
\lambda_R+\frac{(2\sqrt{d_1r_1(1-a_1)}-\vp)^2}{4d_1}< r_1(1-a_1).
\end{align*}
Assume that $\Omega$ satisfies \eqref{cyc} for some $A\in\R$,
$x_0\in\R^N$ and $R>R_0$.
Fix any $R'\in[R_0,R)$ and set
$$
z_0=x_0-(x_0\cdot e)e+(A+1+R')e.
$$
It follows from \eqref{cyc} that
$$
\forall\ s\geq0,\ \
\Omega\supset \bar{B(z_0+se,R')}.
$$

Let $(u,v)$ be the solution of \eqref{clv} in $[0,+\i)\times\bar\Omega$ with initial value $(u_0,v_0)\in\Theta$.
Since $u$ is continuous and positive (from \eqref{0uv1}) on $(0,+\i)\times\bar\Omega$,
one can choose $\eta>0$ small enough so that
$$
\forall \ x\in\bar {B_{R'}},\ \
u(1,x+z_0+(2\sqrt{d_1r_1(1-a_1)}-\vp)e)\geq\eta e^{-\frac{(2\sqrt{d_1r_1(1-a_1)}-\vp)e\cdot x}{2{d_1}}}\psi_{R'}(x)=: w_0(x)
$$
and $w_0\leq1$ in $\bar{B_{R'}}$.
From the choice of $R_0$, even if it means decreasing $\eta$, the function $w_0$ satisfies
$$
d_1\Delta w_0+(2\sqrt{d_1r_1(1-a_1)}-\vp)e\cdot \nabla w_0+r_1w_0(1-a_1-w_0)\geq0
\ \text{ in }B_{R'}.
$$
Let $w$ be the solution of
$$
\bc
w_t=d_1 \Delta w +(2\sqrt{d_1r_1(1-a_1)}-\vp)e\cdot\nabla w+r_1w(1-a_1-w)&\ \text{ in } \ (0,+\i)\times  B_{R'},\\
w=0 &\ \text{ on } \ (0,+\i)\times  \p B_{R'},\\
w(0,\cdot)=w_0&\ \text{ in }\  B_{R'}.
\ec
$$
By the maximum principle, one infers that
$0\leq w\leq1$ in $[0,+\i)\times\bar{B_{R'}}$, and $w$ is nondecreasing in $t$ for all $x\in\bar{B_{R'}}$.
From standard parabolic estimates, $w(t,x)\to w_\i(x)$ uniformly in $x\in\bar{B_{R'}}$ as $t\to+\i$,
where $w_\i$ satisfies the corresponding elliptic equation and $w_\i(x)\geq w_0(x)$ for all $x\in\bar{B_{R'}}$.
Furthermore, the strong elliptic maximum principle implies that
$w_\i>0$ in $B_{R'}$.

For all $t\geq0$ and $x\in\bar\Omega$, define the function $(\tilde u,\tilde v)$ by
$$
\bc
\tilde u(t,x)=u\e(t+1,x+z_0+(2\sqrt{d_1r_1(1-a_1)}-\vp)(t+1)e\r),\\
\tilde v(t,x)=1-v\e(t+1,x+z_0+(2\sqrt{d_1r_1(1-a_1)}-\vp)(t+1)e\r).
\ec
$$
Since $\tilde v\geq0$ in $[0,+\i)\times\bar\Omega$ by $v\leq1$, the function
$\tilde u$
satisfies
\begin{align*}
\tilde u_t
&\geq d_1 \Delta\tilde u +(2\sqrt{d_1r_1(1-a_1)}-\vp)e\cdot\nabla \tilde u+r_1\tilde u(1-a_1-\tilde u)
\end{align*}
for all $t\geq0$ and $x\in\bar{B_{R'}}$.
By \eqref{0uv1}, $\tilde u\geq0$ on $[0,+\i)\times\p{B_{R'}}$.
Then the comparison principle implies that
$\tilde u\geq w$
in $[0,+\i)\times\bar{B_{R'}}$.

Therefore, one has
$$
 \forall\  x\in {B_{R'}},\ \ \liminf_{t\to+\i} (|\tilde u(t,x)|+|\tilde v(t,x)|)\geq w_\i(x)>0.
$$
It implies that $w^*(e,z_0,u_0,v_0)\geq2\sqrt{d_1r_1(1-a_1)}-\vp$ for all $(u_0,v_0)\in\Theta$.
By the choice of $z_0$, there holds
$w^*(e,z,u_0,v_0)\geq2\sqrt{d_1r_1(1-a_1)}-\vp$ for all $(u_0,v_0)\in\Theta$
and $z\in\R^N$
such that $|z-x_0-((z-x_0)\cdot e)e|<R'$.
Since $R'\in[R_0,R)$ is arbitrary, it follows that
$$
w^*(e,z,u_0,v_0)\geq2\sqrt{d_1r_1(1-a_1)}-\vp
$$
for all $(u_0,v_0)\in\Theta$ and $z\in\R^N$
such that $|z-x_0-((z-x_0)\cdot e)e|<R$.
By Proposition \ref{dz}, one also gets that
$$
w^*(e,u_0,v_0)\geq 2\sqrt{d_1r_1(1-a_1)}-\vp
$$
for all $(u_0,v_0)\in\Theta$.
The proof of Theorem \ref{half-cyl} is complete.
\end{pr}

\section{Domains with zero or infinite spreading speeds}\label{S3}
This section is devoted to constructing some special domains for which the spreading speeds might be zero, infinite,
or may depend on position $z$. Namely, we prove the second assertion of part (a) in Proposition \ref{dz} and Theorem \ref{00}.

\begin{pr}[Proof of Proposition \ref{dz}]
We first construct a domain $\Omega$ such that its complement $\R^N\backslash\Omega$ has the shape of an infinite comb with larger and larger teeth.
  Up to translations and rotations, we can assume without loss of generality that
  $e=(1,0,\cdots,0)$ and $z=(0,2,0,\cdots,0)$.
  Pick a sequence of positive real numbers $(a_n)_{n\in\mathbb N}$ such that
  $a_n/n\to+\i$ as $n\to+\i$.
  Define two subsets $\Gamma$ and $ \widehat\Omega$ of $\R^2$ (see Figure \ref{fig5}) as
  $$
  \Gamma=\{(x_1,0):x_1\geq0\}\cup\bigcup_{n\in\mathbb N}(\{n\}\times[0,a_n])\subset\R^2
  \ \text{ and }\
  \widehat\Omega:=\e\{x\in\mathbb R^2:|x-y|<\frac{1}{3},\ \forall\  y \in\Gamma\r\}.
  $$
  \begin{figure}[htbp]
  \centering
  \includegraphics[width=6cm]{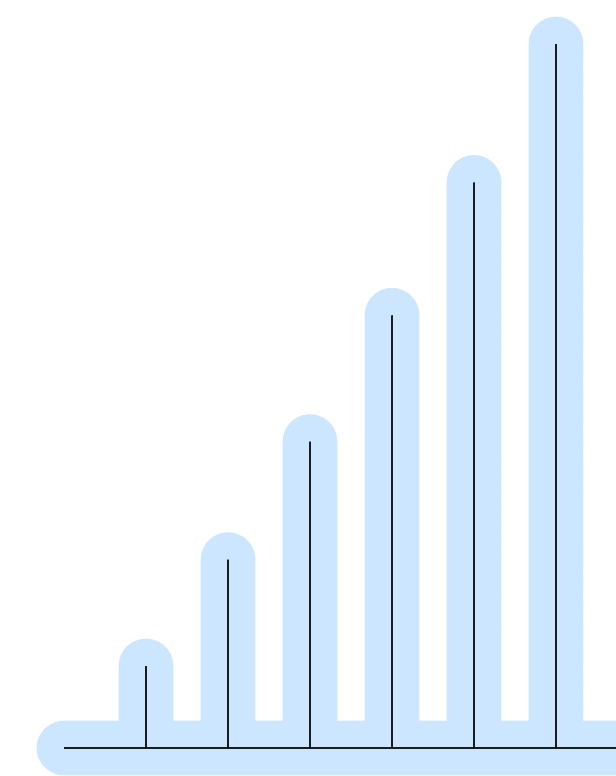}
  \caption{the  line represents $\Gamma$ and the shaded area represents $\widehat\Omega$.}\label{fig5}
\end{figure}Let $\tilde\Omega$ be any open subset of $\R^2$ such that
  $
  \Gamma\subset\tilde\Omega\subset\widehat\Omega
  $
  and such that $\Omega_2:=\mathbb R^2\backslash \bar{\tilde\Omega}$ is connected, locally $C^2$
  and satisfies the extension property.
  Set $\Omega=\Omega_2$ if $N=2$
  and $\Omega=\Omega_2\times\R^{N-2}$ if $N\geq3$.
  It is evident that $\Omega$ is strongly unbounded in direction $e$ and that $\Omega$ does not satisfy
  the assumptions of Theorem \ref{idz}.

 Let $\gamma>0$ be any fixed constant and let $(u_0,v_0)\in\Theta$.
  Remember that $z=(0,2,0,\cdots,0)$. From the construction of $\Omega$, there holds
  $$
  \forall\ s\geq0,\ \
  \bar{B(z+se,1)}\cap\bar\Omega\neq\emptyset.
  $$
  Let $C_0>4$ be given.
 Define $w(0,\cdot)=\max(u_0,1-v_0)$ in $\bar\Omega$.
  Since $d_1=d_2$ and $a_1a_2\leq1$, by \cite[Proposition 2.5]{bhn} and similar arguments as    the proof of Theorem  \ref{upd},
  there exist positive constants $C$ and $\dl$ such that
  \begin{align}\label{he1}
  (0,0)&\leq( u(t,x),1-v(t,x))\notag\\
  &\leq C e^{r_1(1-a_1)t}\|w(0,\cdot)\|_{L^\i(\bar\Omega)}(1+(d_1\dl t)^{-N/2})\int_{{\rm supp} (w_0)}e^{-\frac{d_\Omega(x,y)^2}{C_0d_1t}}dy
  \times(1,a_2)
  \end{align}
  for all $t\geq 0$ and $x\in\bar \Omega$.
  Since ${\rm supp}(u_0)$ and  ${\rm supp}(1-v_0)$ are compact,
  it can be inferred from the definition of $\Omega$ that
  $$
  \inf_{y\in{\rm supp}(w(0,\cdot)),\ s\geq\gamma t,\ x\in\bar{B(z+se,1)}\cap\bar\Omega}
  \frac{d_\Omega(x,y)}{t}\to+\i\ \text{ as }\ t\to+\i.
  $$
Then there exists $t_0>0$ such that
   $$
  (0,0)\leq (u(t,x),1-v(t,x))\leq C e^{-r_1(1-a_1) t}\|w(0,\cdot)\|_{L^\i(\bar\Omega)}(1+(d_1\dl t)^{-N/2})|{\rm supp} (w(0,\cdot))|\times\left( 1,a_2\right)
  $$
  for all $t\geq t_0$, $s\geq \gamma t$ and $x\in\bar{B(z+se,1)}\cap\bar\Omega$,
  where $|{\rm supp} (w(0,\cdot))|$ is the Lebesgue measure of ${\rm supp} (w(0,\cdot))$.
Letting $t\to+\i$, then
$$
\sup_{s\geq\gamma t,\ x\in\bar{B(z+se,1)}\cap\bar\Omega}(|u(t,x)|+|1-v(t,x)|)\to0\ \text{ as }t\to+\i.
$$

Since $\gamma>0$ is arbitrary, thanks to \eqref{w}, there holds
$ w^*(e,z,u_0,v_0)=0
$
for all $(u_0,v_0)\in\Theta$.
On the other hand, the same arguments lead to
$w^*(e,z',u_0,v_0)=0
$
for all $(u_0,v_0)\in\Theta$ and $z'=(z_1',\cdots,z_N')\in\R^N$ such that $z_2'>1/2$.
By Corollary \ref{coro2}, one concludes that
$$
0=w^*(e,z',u_0,v_0)<w^*(e,y,u_0,v_0)=2\sqrt{d_1r_1(1-a_1)}
$$
for all $(u_0,v_0)\in\Theta$,  $z'=(z_1',\cdots,z_N')\in\R^N$ such that $z_2'>1/2$ and $y=(y_1,\cdots,y_N)\in\R^N$ such that $y_2\leq-1/{2}$.
This completes the proof.
\end{pr}

\begin{pr}[Proof of Theorem \ref{00}]
The proof is divided into two steps.

{\bf Step 1: proof of part (a).} We first construct a domain which has the shape of a spiral.
Let $\Omega\subset \R^N$ be a locally $C^2$ connected domain satisfying the extension property
and such that
$\Omega\backslash \bar {B_{2\pi}}=\{x\in\R^2:|x-y|<1,\ y\in\Gamma\}$,
where $\Gamma=\{(t\cos t, t\sin t):t\geq0\}$ (see Figure \ref{fig6}).
Then $\Omega$ is strongly unbounded in every direction $e\in\mathbb S^1$.
\begin{figure}[htbp]
  \centering
  \includegraphics[width=7.5cm]{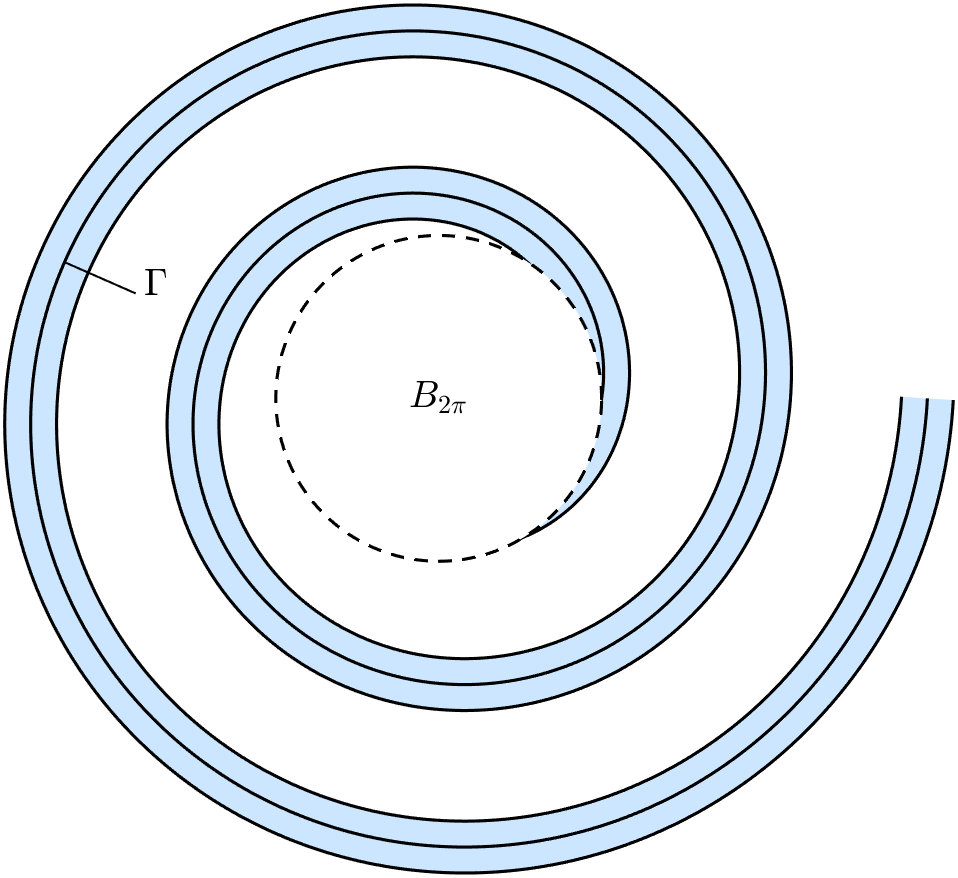}
  \caption{the shaded area is $\Omega\backslash \bar {B_{2\pi}}$, which has the shape of a spiral.}\label{fig6}
\end{figure}

Let $(u,v)$ be the solution of \eqref{clv} in $[0,+\i)\times\bar\Omega$ with any given initial value $(u_0,v_0)\in\Theta$.
Let $e\in\mathbb S^1$ and $C_0>4$ be given.
Clearly, there exists $R>0$ such that $\bar\Omega\cap\bar{B(se,R)}\neq\emptyset$
for all $s\geq0$.
Define $w(0,\cdot)=\max(u_0,1-v_0)$ in $\bar\Omega$.
  Since $d_1=d_2$ and $a_1a_2\leq1$, from  \cite[Proposition 2.5]{bhn} and similar arguments as the proof of Theorem  \ref{upd},
  there exist positive constants $C$ and $\dl$ such that
 \eqref{he1}
 holds for all $t>0$ and $x\in\bar\Omega$.
  Let $\gamma>0$ and $A\geq R$ be any fixed constants.
  Then for all $s\geq0$ and $t>0$, there holds
  \begin{align*}
  (0,0)&\leq \sup_{x\in\bar{B(se,A)}\cap\bar\Omega}(u(t,x),1-v(t,x))\\
  &\leq
  C e^{r_1(1-a_1)t}\|w(0,\cdot)\|_{L^\i(\bar\Omega)}\times(1+(d_1\dl t)^{-N/2})\int_{{\rm supp} (w_0)}e^{-\frac{r_{y,s}^2}{C_0d_1t}}dy\times(1,a_2),
  \end{align*}
  where
  $r_{y,s}=\min_{x\in\bar{B(se,A)}\cap\bar\Omega}d_\Omega(x,y)$.
  Since ${\rm supp}(u_0)$ and ${\rm supp}(1-v_0)$ are compact,
  it derives from the definition of $\Omega$ that
  there exist $\eta>0$ and $t_0>0$ such that
  $$
  \forall\ t\geq t_0,\
  \forall\ s\geq\gamma t,\
  \forall\ y\in {\rm supp}(w_0),
  \ \ r_{y,s}\geq \eta t^2.
  $$
  It implies that
   \begin{align*}
\sup_{s\geq\gamma t,\ x\in\bar{B(se,A)}\cap\bar\Omega}(|u(t,x)|+|1-v(t,x)|)\to0
\ \text{ as }t\to+\i.
  \end{align*}

As a conclusion,
$w^*(e,z,u_0,v_0)=w^*(e,u_0,v_0)=0$
for all $e\in\mathbb S^1$, $z\in\R^2$ and $(u_0,v_0)\in\Theta$.

{\bf Step 2: proof of part (b).}
Let the dimension $N\geq 2$ be fixed. Call $(x_1,x')$ the coordinates in $\R^N$, where $x'=(x_2,\cdots,x_N)$.
Denote $r'=(x_2^2+\cdots+x_N^2)^{1/2}$.
We now turn to construct a domain $\Omega$ with the shape of an infinite cusp.
For all $s\in\R$, define the function $h(s)=e^{-e^s+s}$.
Denote the sets
$$
\tilde\Omega=\{(x_1,x')\in\R^N:x_1>A,\ 0\leq r'\leq h(x_1)\}
$$
and
$$
\widehat\Omega=\tilde\Omega\cup\{(x_1,x')\in\R^N:A-1\leq x_1\leq A,\ 0\leq r'<1\}
$$
for some $A>0$ (see Figure \ref{fig7}).
\begin{figure}[htbp]
  \centering
  \includegraphics[width=14cm]{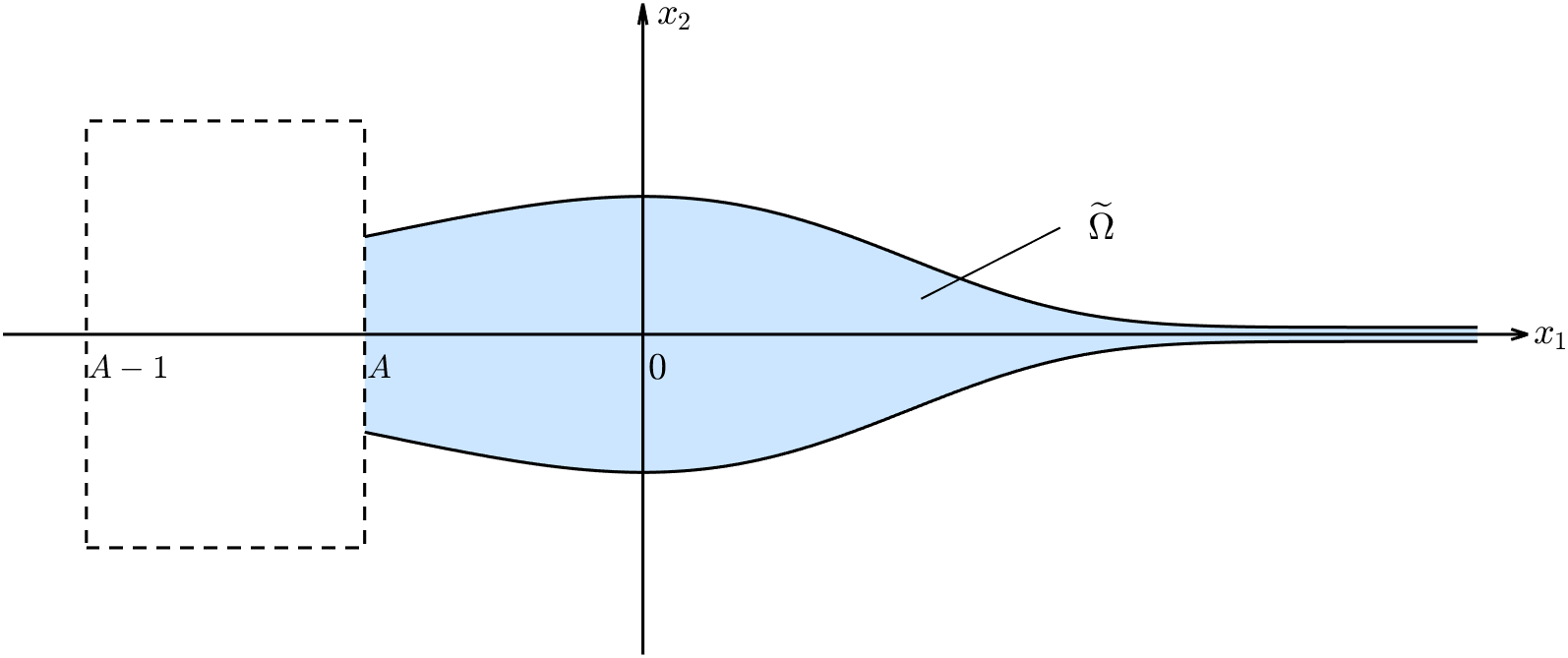}
  \caption{the subsets $\tilde\Omega$ and $\widehat\Omega$ of $\R^2$, the shaded area $\tilde\Omega$ is an infinite cusp.}\label{fig7}
\end{figure}
Let $\Omega$ be an open connected and locally $C^2$ domain such that
$
\tilde\Omega\subset\Omega\subset\widehat\Omega
$.
By increasing $A>0$ if necessary, according to \cite[Lemmas 4.1 and 4.2]{bhn},
there exists $t_0>0$ such that
\begin{align}\label{heatk2}
\inf_{t\geq T,\ y\in K,\ z\in\bar \Omega}p(d_1t,y,z)>0
\end{align}
for all compact subset $K\subset\bar\Omega$, where $p(t,y,z)$ is the heat kernel in $\Omega$
with Neumann boundary condition on $\p\Omega$.
It is easy to see that $\Omega$ is strongly unbounded in the direction $e_1:=(1,0,\cdots,0)$
but $\Omega$ do not satisfy the extension property defined in Section \ref{fund}.

Let $(u,v)$ be the solution of \eqref{clv} with any given initial value  $(u_0,v_0)\in\Theta$.
Call $u_{\rm KPP}$ be the solution of
$$
\bc
(u_{\rm KPP})_t=d_1\Delta u_{\rm KPP}+r_1u_{\rm KPP}(1-a_1-u_{\rm KPP})&\text{ in }\ (0,+\i)\times \Omega,\\
(u_{\rm KPP})_\nu=0&\text{ on }\ (0,+\i)\times\p\Omega,\\
u_{\rm KPP}(0,\cdot)=\min(u_0,1-a_1)&\text{ in }\ \Omega.
\ec
$$
The maximum principle yields that
$0\leq u_{\rm{KPP}}\leq1-a_1$  in $[0,+\i)\times\bar\Omega$.
Since $0\leq v\leq 1$ in $[0,+\i)\times\bar\Omega$ by \eqref{0uv1}, it follows from the comparison principle that
$$
\forall\ t\geq0,\ \forall \ x\in\bar\Omega,\ \ \ u(t,x)\geq u_{\rm KPP}(t,x).
$$
Let $w$ be the solution of  equation $w_t=d_1\Delta w$ with the boundary condition $w_\nu=0$ on $\p\Omega$
and with the initial value   $w(0,\cdot)=u_{\rm KPP}(0,\cdot)$ in $\bar\Omega$.
Applying the comparison principle again,
one has
$$
\forall\ t\geq0,\ \forall \ x\in\bar\Omega,\ \ \
u(t,x)\geq u_{\rm KPP}(t,x)\geq w(t,x).
$$
Note that ${\rm supp} (u_{\rm KPP}(0,\cdot))={\rm supp} (u_0)$.
Then
$$
\forall\ t\geq0,\ \forall \ x\in\bar\Omega,\ \ \
u(t,x)\geq w(t,x)=\int_{{\rm supp }(u_0)}p(d_1t,x,y)u_{\rm KPP}(0,y)dy.
$$
Since ${\rm supp }(u_0)$ is compact,
it follows from \eqref{heatk2} that
there exists $\dl>0$ such that
$p(d_1t,x,y)\geq\dl$ for all $t\geq t_0$, $x\in\bar\Omega$ and $y\in{\rm supp }(u_0)$.
Thus, one has
$$
\forall \ t\geq t_0,\ \forall\ x\in\bar\Omega,\ \ u(t,x)\geq \vp:=\dl \int_{{\rm supp }(u_0)}u_{\rm KPP}(0,y)dy>0.
$$

Let $(\omega_1(t),\omega_2(t))$ be the solution of
 $$
\bc
\dot \omega_1(t)=r_1\omega_1(t)(1-a_1-\omega_1(t)+a_1\omega_2(t)), &t> t_0,\\
\dot \omega_2(t)=r_2(1-\omega_2(t))(a_2\omega_1(t)-\omega_2(t)), &t> t_0,\\
(\omega_1(t_0),\omega_2(t_0))=(\vp,0).
\ec
$$
By the comparison principle, one has
$(u(t,x),1-v(t,x))\geq(\omega_1(t),\omega_2(t))$ for all $t\geq t_0$ and $x\in\bar\Omega$.
Since $(\omega_1(t),\omega_2(t))\to(1,1)$ as $t\to+\i$ (similar to \eqref{ode}) and
 $(u,1-v)\leq(1,1)$ in $[0,+\i)\times\bar\Omega$,
 there holds
$$
(u,v)(t,x)\to(1,0)\ \ \text{ uniformly in }x\in\bar\Omega\text{ as }t\to+\i.
$$
This also  implies that $w_*(e_1,u_0,v_0)$ and $w_*(e_1,z,u_0,v_0)$ are well-defined.

By virtue of  Definitions \ref{s1}-\ref{ss2},
one has
$$w_*(e_1,u_0,v_0)=w_*(e_1,z,u_0,v_0)=w^*(e_1,z,u_0,v_0)=w^*(e_1,u_0,v_0)=+\i$$
for all $(u_0,v_0)\in\Theta$ and $z\in\R^N$.
The proof is complete.
\end{pr}

\section*{Acknowledgments}
The authors sincerely appreciate the valuable and constructive advice provided by professor
Fran\c{c}ois Hamel from Aix Marseille universit\'{e} on this work.
The first author would like to give her sincere thanks to China
Scholarship Council for a 18-month visit of Aix Marseille universit\'{e}.
Her work was also  partially supported by NSF
of China (11971128) and
by the Heilongjiang Provincial Natural Science
Foundation of China (LH2020A003).
The second author's work was partially
supported by NSF of China (12171120).

\section*{Date availability statements}
We do not analyse or generate any datasets, because our work proceeds within a theoretical and mathematical approach.


\begin{thebibliography}{99}
\footnotesize{
\bibitem{AX} M. Alfaro, D. Xiao,
Lotka-Volterra competition-diffusion system: the critical competition case,
{\it Comm. Partial Differential Equations} {\bf48} (2023), 182-208.

\bibitem{AO}A. Alhasanat, C. Ou,
Minimal-speed selection of traveling waves to the Lotka-Volterra competition model,
\textit{J. Differential Equations} \textbf{266} (2019), 7357-7378.

\bibitem{Al}A. Alhasanat, C. Ou, On a conjecture raised by Yuzo Hosono, \textit{J. Dynam. Differential Equations} \textbf{31} (2019), 287-304.


\bibitem{BH1} H. Berestycki, F. Hamel, Generalized traveling waves for reaction-diffusion equations, in Perspectives in Nonlinear Partial Differential Equations. In honor of H. Brezis, American Mathematical Society, Providence, RI, 2007,  101-123.

\bibitem{BH2}H. Berestycki, F. Hamel, Generalized transition waves and their properties, \textit{Commun. Pure Appl. Math.} \textbf{65} (2012), 592-648.




\bibitem{bhn1}H. Berestycki, F. Hamel, N. Nadirashvili,
The speed of propagation for KPP type problems. \uppercase\expandafter{\romannumeral1}. Periodic framework,
\textit{J. Eur. Math. Soc.} \textbf{7} (2005), 173-213.

\bibitem{bhn} H. Berestycki, F. Hamel, N. Nadirashvili,
The speed of propagation for KPP type problems. \uppercase\expandafter{\romannumeral2}. General domains, \textit{J. Amer. Math. Soc.} \textbf{23} (2010), 1-34.

\bibitem{C}C. Carr\`{e}re,
Spreading speeds for a two-species competition-diffusion system,
{\it J. Differential Equations} {\bf264} (2018), 2133-2156.

\bibitem{CG} C. Conley, R. Gardner, An application of the generalized Morse index to travelling wave solutions of a competitive reaction-diffusion model,
{\it Indiana Univ. Math. J.} {\bf33} (1984), 319-343.

\bibitem{D} E. B. Davies, Heat kernels and spectral theory,
{Cambridge University Press}, Cambridge, 1989.

\bibitem{FT}P. C. Fife, M. M. Tang,
Comparison principles for reaction-diffusion systems: irregular comparison functions and applications to questions of stability and speed of propagation of dist,
{\it J. Differential Equations} {\bf40} (1981), 168-185.


\bibitem{Ga} R. A. Gardner,
Existence and stability of travelling wave solutions of competition models: a degree theoretic approach,
{\it J. Differential Equations} {\bf44} (1982), 343-364.

\bibitem{GLam}L. Girardin, K.-Y. Lam,
Invasion of open space by two competitors: spreading properties of monostable two-species competition-diffusion systems,
{\it Proc. Lond. Math. Soc.} {\bf119} (2019), 1279-1335.

\bibitem{G}M. Gruber,
Harnack inequalities for solutions of general second order parabolic equations and estimates of their H\"{o}lder constants,
{\it Math. Z.} {\bf 185} (1984), 23-43.

\bibitem{GL}J.-S. Guo, X. Liang,
The minimal speed of traveling fronts for the Lotka-Volterra competition system,
\textit{J. Dynam. Differential Equations} \textbf{23} (2011),  353-363.

\bibitem{Ho}Y. Hosono, The minimal speed of traveling fronts for diffusive Lotka-Volterra competition model,
\textit{Bull. Math. Biol.} \textbf{60} (1998), 435-448.

\bibitem{Hu}W. Huang, Problem on minimum wave speed for Lotka-Volterra reaction–diffusion competition model, \textit{J. Dynam.
Differential Equations} \textbf{22} (2010), 285-297.

\bibitem{HH}W. Huang, M. Han, Non-linear determinacy of minimum wave speed for Lotka-Volterra competition model, \textit{J. Differential
 Equations} \textbf{251} (2011), 1549-1561.

 \bibitem{HO}Z. Huang, C. Ou, Speed determinacy of traveling waves to a stream-population model with Allee effect,
{\it SIAM J. Appl. Math.} {\bf80} (2020), 1820-1840.



\bibitem{K1}
Y. Kan-on, Parameter dependence of propagation speed of travelling waves for competition-diffusion equations,
{\it SIAM J. Math. Anal.} {\bf26} (1995), 340-363.

\bibitem{K}Y. Kan-on,
Fisher wave fronts for the Lotka-Volterra competition model with diffusion,
\textit{Nonlinear Anal.} \textbf{28} (1997), 145-164.

\bibitem{KF}
Y. Kan-on, Q. Fang, Stability of monotone travelling waves for competition-diffusion equations,
{\it Japan J. Indust. Appl. Math.} {\bf13} (1996), 343-349.

\bibitem{llw} M. A. Lewis, B. Li, H. F. Weinberger,
Spreading speed and linear determinacy for two-species competition models,
\textit{J. Math. Biol.} \textbf{45} (2002), 219-233.





\bibitem{lib}G. M. Lieberman,
Second order parabolic differential equations,
World Scientific Publishing Co., Inc., River Edge, NJ, 1996.

\bibitem{LL}G. Lin, W.-T. Li,
Asymptotic spreading of competition diffusion systems: the role of interspecific competitions,
{\it European J. Appl. Math.} {\bf23} (2012),  669-689.

\bibitem{LLL1} Q. Liu, S. Liu, k.-Y. Lam, Asymptotic spreading of interacting species with multiple fronts \uppercase\expandafter{\romannumeral1}: a geometric optics approach,
{\it Discrete Contin. Dyn. Syst.} {\bf40} (2020), 3683-3714.

\bibitem{LLL2}S. Liu, Q. Liu, k.-Y. Lam,  Asymptotic spreading of interacting species with multiple fronts \uppercase\expandafter{\romannumeral2}: Exponentially decaying initial data,
{\it J. Differential Equations} {\bf303} (2021), 407-455.

\bibitem{MO}M. Ma, C. Ou,
Linear and nonlinear speed selection for mono-stable wave propagations,
{\it SIAM J. Math. Anal.} {\bf51} (2019), 321-345.

\bibitem{MS} R. H. Martin, H. L. Smith,
Abstract functional-differential equations and reaction-diffusion systems,
{\it Trans. Amer. Math. Soc.} {\bf321} (1990), 1-44.



\bibitem{PWZ}R. Peng, C.-H. Wu, M. Zhou,
Sharp estimates for the spreading speeds of the Lotka-Volterra diffusion system with strong competition,
{\it Ann. Inst. H. Poincar\'{e} C Anal. Non Lin\'{e}aire} {\bf38} (2021), 507-547.


\bibitem{PM} M. H. Protter, H. F. Weinberger,
Maximum principles in differential equations,
Springer-Verlag, New York, 1984.

\bibitem{RH} L. Roques, Y. Hosono, O. Bonnefon, T. Boivin,
The effect of competition on the neutral intraspecific diversity of invasive species,
\textit{J. Math. Biol.} \textbf{71} (2015), 465-489.

\bibitem{shen} W. Shen, Traveling waves in diffusive random media, \textit{J. Dynam.
Differential Equations} \textbf{16} (2004), 1011-1060.

\bibitem{st}E. M. Stein, Singular integrals and differentiability properties of functions,
{\it Princeton Math. Ser.} {\bf30}
Princeton University Press, Princeton, NJ, 1970.

\bibitem{TF}M. M. Tang, P. C. Fife, Propagating fronts for competing species equations with diffusion,
{\it Arch. Rational Mech. Anal.} {\bf73} (1980), 69-77.

\bibitem{VVV} A. I. Volpert, V. A. Volpert, V. A. Volpert,
Traveling wave solutions of parabolic systems,
\textit{Transl. Math. Monogr.} \textbf{140} American Mathematical Society, Providence, 1994.
    }

\bibitem{wll}H. F. Weinberger, M. A. Lewis, B. Li, Analysis of linear determinacy for spread in cooperative models,
\textit{J. Math. Biol.} \textbf{45} (2002), 183-218.
\end{thebibliography}
\end{document}